\documentclass[11pt]{article}
\usepackage[top=35mm, bottom=40mm, left=25mm, right=25mm]{geometry}
\usepackage{mathptmx}      % use Times fonts if available on your TeX system
\usepackage{amsmath}
\usepackage{amssymb}

\usepackage{eucal}
\usepackage{paralist}
\usepackage{graphics}                      %% tested in 11 point.\textbf{}
\usepackage{graphicx}
\usepackage{epsfig} %For pictures: screened artwork should be set up with an 85 or 100 line screen

\usepackage{authblk}

\usepackage{hyperref,setspace}
\usepackage{enumerate,verbatim,natbib}

\usepackage[title,titletoc]{appendix}

\usepackage{setspace}
\usepackage{indentfirst}
\allowdisplaybreaks[1]%3

\newtheorem{theorem}{Theorem}[section]
\newtheorem{proposition}{Proposition}[section]

\newtheorem{lemma}{Lemma}[section]
\usepackage{color}
\usepackage[misc,geometry]{ifsym}
\usepackage{natbib}

\newcommand{\var}{\varepsilon}

\newcommand{\p}{\partial}

\title{Modeling the distribution of insulin in pancreas\thanks{JY was supported  by National Natural Science Foundation of China 61573016  when the work was performed. JJ is supported by the Juvenile Diabetes Research Foundation. JL was supported in part by NIH Grant R01-DE019243 and DOE Grant DE-EM000197 when the work was performed.}}
\author[a]{Changbing Hu}
\author[b]{Junyuan Yang}
\author[c]{James D. Johnson}
\author[a,\thanks{Corrsponding author: jiaxu.li@louisville.edu}]{Jiaxu Li}

\affil[a]{Department of Mathematics, University of Louisville, Louisville KY 40292, USA}
\affil[b]{Complex Systems Research Center, Shanxi University, Taiyuan 030006, P.R. China}
\affil[c]{Department of Cellular and Physiological Sciences and Department of Surgery, University of British Columbia, Vancouver, BC, Canada, V6T 1Z3}

\begin{document}
	
\maketitle

%**********************************************************************
\begin{abstract}
Maintenance of adequate physical and functional pancreatic $\beta$-cell mass is critical for the prevention of diabetes mellitus or to postpone its onset. It is well established that insulin potently activates mitogenic and anti-apoptotic signaling cascades in cultured $\beta$-cells. Furthermore, loss of $\beta$-cell insulin receptors is sufficient to induce type 2 diabetes in mice. However, it remains unclear whether the {\em in vitro} effect in human islets and the {\em in vivo} effects in mice can be applied to human physiology. The major obstacle to a complete understanding of the effects of insulin's feedback in human pancreas is the absence of technology to measure the concentrations of insulin inside of pancreas. To contextualize recent {\em in vitro} data, it is essential to know the local concentration and distribution of insulin in pancreas. To this end, we continue to estimate the local insulin concentration within pancreas. In this paper, we investigate the distribution of insulin concentration along the pancreatic vein through a novel mathematical modeling approach using existing physiological data and islet imaging data, in contrast to our previous work focusing on the insulin level within an islet. Our studies suggest that, %{\textcolor{red}
in response to an increase in glucose, the insulin concentration along the pancreatic vein increases nearly linearly in the fashion of increasing quicker in tail area but slower in head area depending of the initial distribution.%}
Our studies also reveal that the distribution of islets in pancreas might be significant for the steady state insulin concentrations in pancreas vein. A widely believed statement, insulin concentration in pancreas is much higher than periphery, could be only true for some cases of islet distribution. We also find that the factor of small diffusion with blood is negligible since the convection of blood flux dominates. Our work might also shed a light to the future design of islet implantation surgeries in regarding to whether and how to arrange the distribution of implanted islets. To the best of our knowledge, our model is the first attempt to estimate the insulin distribution in pancreatic vein. Our model is simple, robust and thus can be easily adopted to study more sophistic cases. 
\end{abstract}
%**********************************************************************

%\MSC[2008] 34K20 \sep 92C50 \sep 92D25
%\keywords

\noindent{{\bf \emph{Keywords:}} { Pancreas, \and Islet, \and Insulin distribution, \and splenic vein, \and Initial-boundary value problem, \and Singular perturbation, \and Steady state solution }

%, the main secretory product of $\beta$-cells,
%, and the dose-response relationships between insulin and $\beta$-cell survival

%================================================================
%---------------------------------------------------
\section{Introduction} \label{Sec:intro}
%---------------------------------------------------

Diabetes mellitus continues to exact a devastating toll on society due to its life-threatening complications. Diabetes mellitus is a leading cause of heart disease, kidney failure, blindness and amputations, as well as other pathologies (\cite{derouich2002effect}). At least 5\% of the worlds population has diabetes, and up to 30\% of the population may be at risk. Despite decades of study, the factors controlling the initiation and progression of diabetes still remain to be fully elucidated. Without a thorough understanding of the glucose homeostasis system and its dysfunction in diabetes, we will continue to struggle to develop new approaches to detect, prevent and delay the onset of diabetes. Thus, it is important to integrate important but reductionist experimental findings into comprehensive models. Therefore, there is a pressing need for accurate mathematical models employing the latest experimental findings.

Diabetes mellitus is characterized by either the near total loss of pancreatic $\beta$-cells (type 1 diabetes), or the relative loss of functional $\beta$-cells that accompanies the inability of the body's cells to utilize insulin properly (type 2 and gestational diabetes). The pancreatic $\beta$-cell is the unique cell type that secretes insulin, a hormone needed for most cells to uptake and store glucose. Hormones such as insulin ensure that the plasma glucose concentration remains in the normal range during and between meals. In addition to dysfunctional insulin secretion (\citet{bergman2002evolution, lowell2005mitochondrial}), it is now widely accepted that inadequate compensatory changes in pancreatic $\beta$-cell mass plays a critical role at the onset of type 2 diabetes and gestational diabetes (\citet{ackermann2007molecular, rhodes2005type}). Thus, it becomes important and urgent to better understand the dynamics of $\beta$-cell proliferation and programmed cell death in the regulation of $\beta$-cell mass, and the causes of failure for $\beta$-cell adaptation. These findings will aid in developing effective and efficient therapeutic agents.

The physical and functional $\beta$-cell mass is determined by the number and sizes of all beta-cells, as well as their differentiation state. Numerous experiments indicate that insulin itself has anti-apoptotic and proliferative effects on $\beta$-cells (\citet{alejandro2010acute, beith2008insulin, johnson2006insulin}), %\textcolor{red}{
beta-cell mass 2010.
%}. 
Importantly, it has been shown that insulin has a bell shaped dose-response curve \cite{johnson2008control}.% \textcolor{red}{[XXX also cite Johnson 2002 PNAS]}. 
In order to discover the importance of these effects, it is critical to understand the distribution of insulin concentrations in pancreas. In this paper, we propose a partial differential equation (PDE) model to investigate the insulin distribution in pancreas. We organize this paper as follows. In the next section, we formulate the PDE model based on evidence from physiological studies and consider the existence of  mild solutions of the PDE model in Section \ref{Sec:model}, the existence of the steady state solution and that the steady state solution is stable in Section \ref{sect:main}. In Subsection \ref{subsect:epsilon}
, we show that the convection of the blood flux in pancreas dominates over the small diffusion effect. Numerical simulations are carried out in Section \ref{Sec:sim}. We end this paper with a discussion of the physiological ramifications in Section \ref{Sec:disc}. The details of theoretical proof are enclosed in Appendix.% Section \ref{sect:proof}.

%---------------------------------------------------
%================================================================

%================================================================
\section{Model formulation} \label{Sec:model}
%---------------------------------------------------

Pancreas is an elongated and tapered organ that is located across the back of the abdomen and behind the stomach. The pancreas has been divided into three morphological sections: the head, the body and the tail. Arterial blood flows into the tail and flows out of the head of the pancreas through the hepatic portal vein. After insulin flows out of the pancreas through the hepatic portal vein (the nearest location to the islets at which insulin concentration is known), the liver and kidney clear about 50\% of the secreted insulin before it circulates around the whole body, becoming what we will refer to herein as peripheral insulin. The unused (undegraded) peripheral insulin returns to the pancreas with the blood influx. The half-life if insulin in the blood is ~5 minutes, meaning that most of the insulin in the body is turned over every 10 minutes.

%------------------------------------------
\begin{figure} [!htbp]
	\centerline{\resizebox{4.5in}{2.5in}{\includegraphics{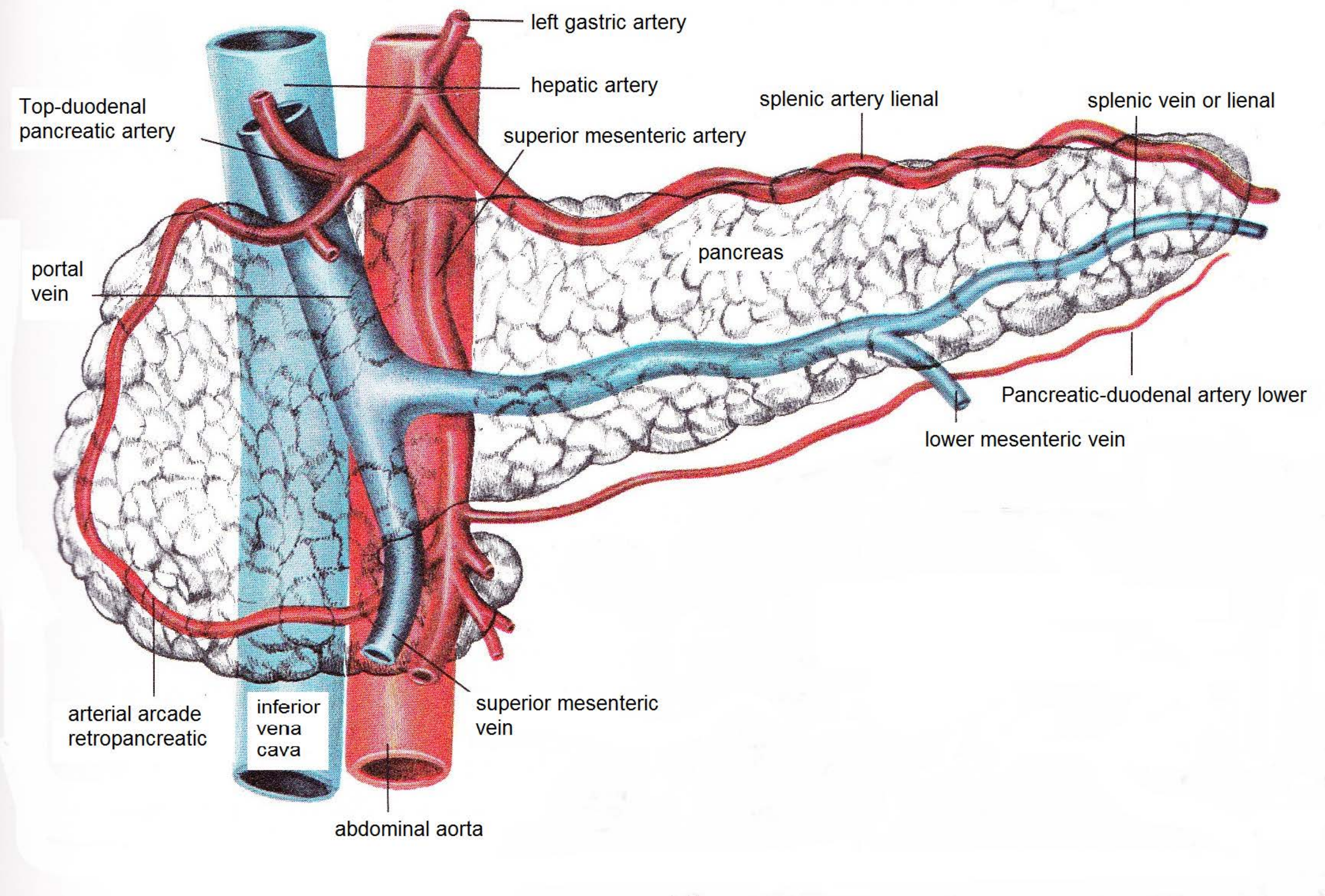}}}
	\caption{\footnotesize \hskip 0.2in Pancreatic blood flow diagram. The pancreatic vein is in blue.}
	\label{PancBloodFlowDiagram}
\end{figure}
%------------------------------------------

Insulin-secreting $\beta$-cells are contained in micro-organs called the islets of Langerhans, spread out within the pancreas and accounting for approximately 1\% of the mass of the entire pancreas. Each islet contains a few hundred to a few thousand $\beta$-cells. There are about one million of islets in an adult human pancreas. Insulin is secreted in response to elevated glucose in pusatile bursts and oscillations. When diabetes develops, either all of or most of the $\beta$-cells become dysfunctional or the compensation of $\beta$-cell mass is not enough to make up for insulin resistance. Normally, dynamic $\beta$-cell mass requirements are balanced by changes in cell number (through proliferation and/or programmed cell death) and/or cell size (through hypertrophy and atrophy). Generally, the acute effects of glucose and insulin on these dynamics are difficult to observe directly in humans, especially with the current state of medical imaging. Particularly, $\beta$-cell proliferation or death cannot be directly imaged non-invasively, due in part to the very low frequency of dividing cells in this largely 'post-mitotic' tissue. Proliferation is dramatically lower in older animal models and is hardly seen after the age of thirty in humans (\cite{perl2010significant}).% \textcolor{red}{[XXX also Cite papers from Kushner on proliferation in older mice and cite our new review article in Diabetes as Szabat et al in press]}. 
Therefore, investigating the factors that contribute to $\beta$-cell proliferation and apoptosis by mathematical modeling is a relevant and effective tool in understanding the cause of improper adaptation in $\beta$-cell mass. To the best of our knowledge, no mathematical model has been developed to study the effects of local insulin concentrations on pancreatic islet survival and proliferation. These effects may be physiologically significant if they increase $\beta$-cell mass by slowing down the death of the cells, and therefore postpone type 2 diabetes progression. Here, we propose a PDE model to study the distribution of insulin in pancreas along the flow of blood. Glucose and peripheral insulin enter the pancreas with the blood flow and are distributed through capillaries. The insulin in the pancreas is a composite of the peripheral insulin with the insulin newly synthesized and released from $\beta$-cells.

%\textcolor{red}{[XXX this is a small paragraph that should be integrated somewhere else]}
%\textcolor{blue}{The minimal need of insulin (basal insulin) is secreted from $\beta$-cells in pulsatile manner every a few minutes. Elevated glucose at mealtime triggers insulin secretion in an oscillatory fashion which can last between 60 and 150 minutes.}

%------------------------------------------
\begin{figure} [!htbp]
	\centerline{\resizebox{4.5in}{2.5in}{\includegraphics{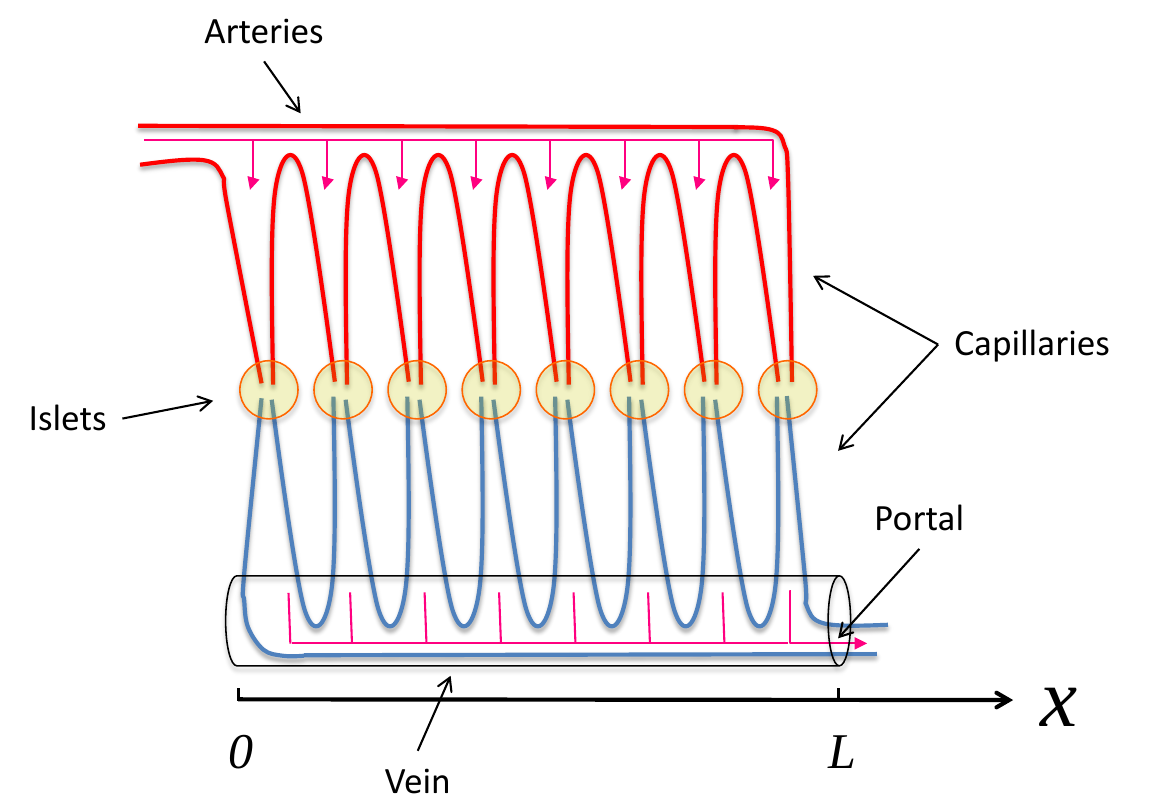}}}
	\caption{\footnotesize \hskip 0.2in Model diagram.}
	\label{ModelDiagram}
\end{figure}
%------------------------------------------

The elongated shape of pancreas and the distribution of many islets as `beads on a string' allow us to assume that a pancreas occupies a one-dimensional space with length $L>0$ (cm) with the tail at the original and the head at the point $L$. (Refer to the clear images in \citet[Fig. 2]{hornblad2011improved}.)  We assume that the blood flow in a small individual organ approximately keeps at a constant rate. Hence, we can assume that the blood flux is at a constant speed $c>0$ from the original $O$ to $L$. Islets, and therefore $\beta$-cells, are distributed along the interval $(0, L)$. Islets are not randomly distributed (\cite{hornblad2011improved}). However, for simplicity, we assume that the distribution is homogeneous along this line. At the same time, glucose and peripheral insulin flow in to pancreas with blood along this line. Inside pancreas, both glucose utilization and insulin degradation are very small relative to the whole body, but they still exist. Since we consider the dynamics of insulin and glucose distribution inside pancreas, both the hepatic glucose production and the insulin-independent glucose utilization, for example, by brain cells, are neglected. (Refer to \citet[Fig. 2]{li2006modeling} for a full systemic description.) The effects resulted from these phenomena are incorporated in the boundary conditions in view of physiological evidences. \cite{nyman2010glucose} studied blood flux in eight single $\beta$-cells by imaging. According to the images obtained by \cite{nyman2010glucose}, the blood flow in islets can be classified in three ways: inner-to-out (diffusion) (4 out 8), top-to-down (convection) (3 out of 8), and irregular (1 out 8) (\citet[Fig. 3]{nyman2010glucose}). According to our calculations, the diffusion, if it exists, is very small. Based on available physiological data, we hereby construct the following general partial differential equation (PDE) model in one dimensional space $[0,L]$, $L>0$, with the pancreas tail at the origin and the head at the $L$:
\begin{equation}
\left\{\begin{array}{l}
\displaystyle -\varepsilon \frac{\partial^2 G}{\partial x^2} + \frac{\partial G}{\partial t} + c\frac{\partial G}{\partial x} =  G_{in} - aGI, \\
\\
\displaystyle -\varepsilon \frac{\partial^2 I}{\partial x^2} + \frac{\partial I}{\partial t} + c\frac{\partial I}{\partial x} =  \frac{\sigma G^2}{b^2 + G^2} - d_iI, \ \ x\in(0,L),\quad t>0,\\
\\
G(L,t)=\alpha_1G(0,t), \qquad I(L,t)=\alpha_2 I(0,t),\qquad t>0,\\
\\
G_x(L,t)=G_x(0,t)=0, \qquad I_x(L,t)=I_x(0,t)=0,\qquad t>0,\\
\\
G(x,0)=\phi(x), \qquad I(x,0)=\psi(x), \qquad x\in(0,L),\\
\end{array}\right.\label{Model_pancreas}
\end{equation}
where $c>0$ is the speed of the blood flow in the pancreas, $\varepsilon>0$ is the constant diffusion rate, $\phi(x)$ and $\psi(x)$ indicates the initial distributions of glucose and insulin in pancreas, respectively, while $G_{in}>0$ is the constant average peripheral inputs of glucose. The term $aGI$ stands for the insulin-dependent glucose uptake with $a>0$ is the insulin sensitivity index, the term $\frac{\sigma(x) G^2}{b^2 + G^2}$ stands for the insulin secretion triggered by glucose, which is dependent of the space location $x\in[0, L]$, $b>0$ is the half saturation of the maximum insulin secretion, and $d_i$ is the constant rate of insulin degradation. %\textcolor{red}{
$\alpha_1,\alpha_2\geq 1$ are proportionality constants to be determined by physiology.% fact, see [?????????]}. 
We naturally assume that $\sigma(\cdot)$ is bounded in $[0, L]$ and denote
\[\bar\sigma:=\sup\limits_{x\in[0,L]}\sigma(x)<+\infty,\quad \mbox{ and } \quad \underline\sigma:=\inf\limits_{x\in[0,L]}\sigma(x)>0.\]
Regarding to the proposed boundary conditions, the authors are not aware of any relevant mathematical analysis in the context of PDEs. Much more are involved to deal with the current boundary condition than the classical periodic boundary condition that is a particular case for $\alpha_1=\alpha_2 = 1$ (refer to Section \ref{sect:main}). We analytically study the Model (\ref{Model_pancreas}) for the case $\sigma(x)$ is a constant $\sigma>0$, and numerically investigate the model for various functional forms $\sigma(x)$.

Regarding to the well-posedness of the system \eqref{Model_pancreas}, define the following functional spaces,
\[X=C([0,L],\mathbb R), \quad \quad Y=X^2,\]
and their positive cones 
\[X_+=C([0,L],\mathbb R_+), \quad \quad Y_+=X^2_+\]
with the following norms, respectively,
\[\|\phi\|_X=\max_{x\in[0,L]}\phi(x), \quad \|\psi\|_Y=\sqrt{\|\psi_G\|^2+\|\psi_I\|^2}, \quad\quad \forall \phi\in X, \quad \psi=(\psi_G,\psi_I)^T\in Y,\]
Now define a liner operator $\mathcal A: Y\rightarrow \mathbb R^2$ defined by 
\begin{align}\label{eqn:A}
\mathcal A\psi=\left(\varepsilon\frac{\partial^2}{\partial x^2}-c\frac{\partial}{\partial x}\right)\psi,\psi\in Y,
\end{align}
with its domain in form of 
\[D(\mathcal A)=\{u\in Y|u(L)=Du(0),u'(L)=u'(0)=0\}\]
where 
\[
  D = \left[\begin{array}{cc} \alpha_1 & 0\\
                                        0 & \alpha_2\end{array}\right] .
\]
It is easy to show that $\mathcal A$ generate a semigroup $S(t)$ that satisfies
$$\mbox{(i) } S(0)\psi(x)=\phi(x),\phi\in Y, \quad  \mbox{(ii) } S(t)S(\tau)\phi=S(t+\tau)\phi, \mbox{ for all } 0\le t, \tau<+\infty, \quad  \mbox{(iii) } \lim_{t\rightarrow0} S(t)\phi=\phi,$$
which implies that $S(t)$ is a $C_0-$semigroup. 

Moreover, it is easy to check that $\mathcal A$ and $S(t)$ have the following properties
\begin{proposition}
Let $\mathcal A$  be defined in \eqref{eqn:A} and $\rho(\mathcal A)$ denote the resolvent set of $\mathcal A$. Then we have
\begin{itemize}
\item [(a)] $\|S(t)\|\leq 1,$ and
%\item $(\omega,\infty)\subset \rho(A)$, and for $\lambda > \omega$,
%\[\|(\lambda -A)^{-n}\|\leq \frac{M}{(\lambda-\omega)^n}\]
\item [(b)] $(0,\infty)\subset \rho(A)$, and for $\lambda > 0$,  
\[\|(\lambda -\mathcal A)^{-n}\|\leq \frac{1}{\lambda^n}.\]
\end{itemize}
\end{proposition}
%-----------------------------------------------------------
Let
\[f(u)=
	\left(
		\begin{array}{c}
			G_{in} - aIG\\
			\\
			\frac{\sigma G^2}{(b^2+G^2)} - d_iI
		\end{array}
	\right), \quad \mbox{ where } \quad u = 
	\left(
		\begin{array}{c}
			G\\
			I
		\end{array}
	\right)\in Y.
\]
It is clear that $f$ is locally Lipschitz. According to \citep[Proposition 4.16]{webb2008population}, we have the following result for the system \eqref{Model_pancreas}. 
%***********************************************************
\begin{theorem} For each $u_0\in Y$, there exists a unique continuous
mild solution $u(t, u_0)\in Y$, for $t>0$, such that
\[u(x,t)=(S(t)u_0)(x)+\int^t_0S(t-\tau) f(u(\tau))\,d\tau.\]
\end{theorem}
%***********************************************************
%================================================================

%================================================================

%================================================================
%---------------------------------------------------
\section{Model analysis\label{sect:main}}
%---------------------------------------------------

For simplicity, we assume that the speed of blood flow in pancreatic is $c=1$. Otherwise, this can be attained by a simple transformation $\bar{L} = L/c$ and thus we need only adjust the length of pancreatic vein in numerical work in Section \ref{Sec:sim}.

%---------------------------------------------------
\subsection{Spacial homogeneous solution\label{subsect:ode}}
%---------------------------------------------------

A spacial homogeneous solution $u_0$ of the system \eqref{Model_pancreas} refers to $u_0$ such that $f(u_0)=0$. That is, $u_0$ is an equilibrium point of the following ordinary differential equation (ODE) system stemming from Model \eqref{Model_pancreas} without spatial differentiation with $\sigma(x)$ assumes a constant $\sigma>0$,
%---------------------------------------------------
\begin{equation}
\left\{\begin{array}{l}
\displaystyle \frac{\partial G}{\partial t} =  G_{in} - aGI, \\
\\
\displaystyle \frac{\partial I}{\partial t}  =  \frac{\sigma G^2}{b^2 + G^2} - d_iI,
\quad t>0,
\end{array}\right.\label{Model:ode}
\end{equation}
%---------------------------------------------------
Based on the graph of the growth function in the right hand side of \eqref{Model:ode}, it is straight forward to check that the ODE system \eqref{Model:ode} has a unique steady state $u^\star=(G^\star, I^\star)$. We first find a bound of $u^{\star}$ and then investigate its stability with respect to the ODE system \eqref{Model:ode}. %(a graph can be inserted here).
To study the stability of $u^star$, we consider the linearized matrix $B\in\mathbb R^{2\times 2}$ of the system
\eqref{Model:ode} near $u^\star= (G^\star, I^\star)$, where
\[B=\left[\begin{array}{cc}
-aI^\star & -aG^\star\\
\\
\frac{2\sigma b^2G^\star}{(b^2+G^{\star 2})^2} & -d_i\end{array}\right].\]
The following lemma gives the existence of equilibria of \eqref{Model:ode} and a priori estimation. 
%---------------------------------------------------
\begin{lemma}\label{lem:equiexistence}
The system \eqref{Model:ode} has a unique equilibrium point $u^\star=(G^\star, I^\star)$ and satisfies
\[\frac{G_{in}d_i}{a\sigma}<G^\star
<\max\left(\frac{2G_{in}d_i}{a\sigma},b\right),
\qquad 0<I^\star <\frac{\sigma}{d_i}.\]
Furthermore $u^\star$ is globally stable in $\mathbb R^2_+$, and there exist positive constants $M$ and $\rho$ such that
\[\|e^{Bx}\|_2\leq Me^{-\rho x},\quad x\geq 0,\]
where $\|\cdot\|_2$ denotes the standard matrix norm, in case of vectors, the norm
$\|\cdot\|_2$ represents the standard Euclidean norm.
\end{lemma}
The proof of Lemma \ref{lem:equiexistence} is given in Appendix \ref{Appen:Lem_EquiExistence}. %\ref{sect:proof}.

%---------------------------------------------------
\subsection{Patial differential equations with $\varepsilon=0$\label{subsect:pde}}
%---------------------------------------------------

In this sebsection, we rewrite the system \eqref{Model_pancreas} without second
order diffusion as the abstract evolution equation for the above system: Let $u=(G,I)^t$,
%---------------------------------------------------
\begin{equation}\label{Model:pancreasII}
\frac{\partial u}{\partial t}+\frac{\partial u}{\partial x}=B(u-u^\star)+F(u).
\end{equation}
%---------------------------------------------------

\[F(u)=\left(\begin{array}{c}
-a(I-I^\star)(G- G^\star)\\
\\
\frac{\sigma b^2(G-G^\star)^2}{(b^2+G^{\star 2})^2}\left[1
-\frac{(G+G^\star)^2}{b^2+G^{2}}\right]\end{array}\right).\]
We turn to the existence of the steady solution $u$ to
\eqref{Model:pancreasII}, i.e.,
the first order ODE with boundary value problem:
%%%%%%%%%%%%%%%%%%%%%%%%%%%%%%%%%%%%%%%%%%%%%%%%%%%%%%%%%%%%%
\begin{equation}\label{steadyEqua}\left\{\begin{aligned}
&\frac{d G(x)}{d x}=G_{in}-aG(x)I(x),
\\
&\frac{d I(x)}{d x}=\frac{\sigma(x)G^2(x)}{b^2+G^2(x)}-d_iI(x),
\\
&G(L)=\alpha_1G(0),I(L)=\alpha_2I(0),\end{aligned}\right.
\end{equation}
Solving equation \eqref{steadyEqua} yields
\begin{align}
G(x)=&G_0e^{-\int_0^xI(s)ds}+G_{in}\int_0^xe^{-\int_\tau^xI(s)ds}d\tau,\label{eqn:G}
\\
I(x)=&I_0e^{-d_ix}+\int_0^xe^{-d_i(x-\tau)}\frac{\sigma (\tau)G^2(\tau)}{b^2+G^2(\tau)}d\tau:=I[G](x).\label{eqn:I}
\end{align}
Combining the boundary conditions, we derive the following compatibility conditions
\begin{align}
G_0=\frac{G_{in}\int_0^Le^{-\int_\tau^LI(s)ds}d\tau}{\alpha_1-e^{\int_0^LI(s)ds}},\label{eqn:G_0}
\\
I_0=\frac{\int_0^Le^{-d_i(x-s)}\frac{\sigma(s)G^2(s)}{b^2+G^2(s)}ds}{\alpha_1-e^{-\int_0^LI(\tau)d\tau}}.\label{eqn:I_0}
\end{align}
Based on the positive operator theory, the following lemma gives the existence of the steady state of \eqref{steadyEqua}. The existence of the steady state is  based on Krasnoselskii Fixed Point Theorem \citep{krasnoselskij1964positive}. 
\begin{lemma}\label{lem:existstate}
Suppose that equations \eqref{eqn:G_0} and \eqref{eqn:I_0} hold. If $\alpha_j>1,(j=1,2)$, then \eqref{steadyEqua} has at least one positive solution.
\end{lemma}
Furthermore, provided that %\textcolor{red}{
$L$ is not long enough.
%}, 
we can obtain the uniqueness of the steady state.
\begin{lemma}\label{lem:steadyuniq}
Suppose that the conditions of Lemma \ref{lem:existstate} hold. If $L$ %\textcolor{red}{
is not long enough,%}, 
then \eqref{steadyEqua} has at most one positive solution.
\end{lemma}
We leave the proofs of Lemma \ref{lem:equiexistence} and \ref{lem:steadyuniq} in Appendix \ref{Appen:Lem_ExistState} and \ref{Appen:Lem_SteadyUniq}

In what follows, we assume $\sigma(x)=\sigma$. By the variation of the costant formula, we readily obtain the steady state satisfiying
\begin{align}
\label{fixedPoint}
u(x)=&u_0(x)+\int^x_0e^{B(x-y)}F(u(y))dy\notag
\\
&+e^{Bx}(D-e^{BL})^{-1}\int^L_0e^{B(L-y)}F(u(y))dy.
\end{align}
Now we are in a positon to focus on the fixed point of \eqref{fixedPoint}.  First we consider the case when the spatial length $L\ll 1$, with the explicit expression on $u_0(x)$ we have
\[u_0(0)\approx 0, \qquad u_0(L)\approx 0,\]
which is not physically relevant. While for large spatial length $L$,
we have
\[u_0(0)\approx D^{-1}u^\star, \qquad u_0(L)\approx u^\star,\]
which fits for the context of anatomy and physiology. So in the following, we always assume that the spatial
length $L$ is sufficiently large
such that
\[\|(D-e^{BL})^{-1}(I-D)\|_2\leq 1-\theta,\]
where $0<\theta<1$.
Hence, we can define a nonlinear opertor $\hat{\mathcal T}: Y\rightarrow Y$ defined by
\begin{align*}
\hat{\mathcal T}[u](x)=&u_0(x)+\int^x_0e^{B(x-y)}F(u(y))dy\notag
\\
&+e^{Bx}(D-e^{BL})^{-1}\int^L_0e^{B(L-y)}F(u(y))dy.
\end{align*}
The existence of the steady state indicates that $\hat{\mathcal T} u=u.$ To proceed such issue, it is not hard to see that the nonlinear term $F$ has the following property
\[\|F\|_X\le k[(u-u^*)^2,\]
where $k=a+\frac{\sigma b^2}{(b^2+G^*)^2}.$  
We introduce a ball domain $B_r(u^\star) \subset C([0,L])^2$ as follows
\[B_r(u^\star)=\{u\in C([0,1])^2:\quad \|u-u^\star\|_\infty\leq r\}.\]
We will try to prove that for small $r$, $T B_r(u^\star)\subset B_r(u^\star)$. We first present some estimates.

%---------------------------------------------------
\begin{lemma}\label{lem:TBound}
With the nonlinear operator $T$ and convex ball set $B_r(u^\star)$,
assume that the equilibrium $u^\star$ satisfies
\[\|u_0-u^\star\|_\infty
\leq \|(D-e^{BL})^{-1}(I-D)u^\star\|_2=\|\xi_0^\star\|_\infty,\]
and
\[\begin{split}&\left\|e^{Bx}(D-e^{BL})^{-1}
\int^L_0e^{B(L-y)}F(u(y))dy+\int^x_0e^{B(x-y)}F(u(y))dy\right\|_\infty\\
&\leq Lk(\|(D-e^{BL})^{-1}\|_2+1)\|u-u^\star\|^2_\infty.
\end{split}\]
\end{lemma}
%---------------------------------------------------
%---------------------------------------------------
%\begin{proof}
{\bf{Proof}}
First we notice that
\[\|u_0-u^\star\|_2\leq \|e^{-\rho x}(D-e^{BL})^{-1}(I-D)u^\star\|_2,\]
which yields the first inequality. For the other two terms in $Tu$, we have
\[\begin{split}&\left\|e^{Bx}(D-e^{BL})^{-1}
\int^L_0e^{B(L-y)}F(u(y))dy+\int^x_0e^{B(x-y)}F(u(y))dy\right\|\\
&\leq \|(D-e^{BL})^{-1}\|Lk\|u-u^\star\|_\infty^2+Lk\|u-u^\star\|_\infty^2\\
&\leq Lk(\|(D-e^{BL})^{-1}\|_2+1)\|u-u^\star\|^2_\infty.
\end{split}\]
$\hfill\square$
%\end{proof}
%---------------------------------------------------

%---------------------------------------------------
\begin{theorem}\label{steadySol}Assume that
\[4\|\xi_0^\star\|_\infty Lk(\|(D-e^{BL})^{-1}\|_2+1)\leq 1.\]
Let
\[r=\frac{1-\sqrt{1-4\|\xi_0^\star\|Lk(\|(D-e^{BL})^{-1}\|_2+1)}}{2}.\]
For any $u\in B_r(u^\star)$
\[\|Tu-u^\star\|_\infty\leq \|(D-e^{BL})^{-1}(I-D)u^\star\|_2
+Lk(\|(D-e^{BL})^{-1}\|_2+1)r^2\leq r,\]
then there exist a steady solution denoted by
$\overline{u}(x)$ to the evolution equation \eqref{Model_pancreas}, which satisfies the bound
\[\|\overline{u}-u^\star\|_{\infty}\leq r.\]
\end{theorem}
%---------------------------------------------------
%---------------------------------------------------
%\begin{proof}
{\bf{Proof}}
The existence proof is straightforward by using the classical contraction
mapping theorem with the estimates in Lemma \ref{lem:TBound}. $\Box$
%\end{proof}
%***************************************************

Now, we focus on the stability of \eqref{Model:pancreasII} under the case $\sigma(x)=\sigma.$  First we compute
the linearization around $\overline{u}$.
Let $u(t, x)=\overline{u}(x)+v(t,x)$
\begin{equation}\frac{\p v}{\p t}=-\frac{\p v}{\p x} +\overline{B}(x)v,
\label{lin}\end{equation}
where
\[\overline{B}=\left[\begin{array}{cc}
-a \overline{I} & -a \overline{G}\\
\frac{2\sigma b^2\overline{G}}{(b^2+\overline{G}^2)^2} & -d_i\end{array}\right].\]
The stability will be determined by the spectral
\begin{equation}
-\frac{\p v}{\p x} +\overline{B}(x)v=\lambda v,
\label{spectral}\end{equation}
with the boundary condition
\[v(L)=Dv(0).\]
The general solution can be written as
\[v(x)=e^{\int^x_0(\overline{B}(y)-\lambda)dy} v(0),\]
by the boundary condition, we obtain
\[Dv(0)= e^{\int^L_0(\overline{B}(y)-\lambda)dy} v(0).\]
In order for $\lambda$ to be an eigenvalue, the determinant must vanish, i.e.,
\[\mbox{det}|D-e^{\int^L_0(\overline{B}(y)-\lambda)dy}|=0.\]
Let
\begin{equation}e^{\int_0^L\overline{B}(y)dy}=
\left[\begin{array}{cc}
b_{11} & b_{12}\\
b_{21} & b_{22}\end{array}\right].\label{bbar}\end{equation}
We obtain the equation on $e^{-\lambda L}$:
\[\alpha_1\alpha_2-e^{-L\lambda}(\alpha_1b_{22}+\alpha_2b_{11})+
e^{-2\lambda L}(b_{11}b_{22}-b_{12}b_{21})=0.\]
The the eigenvalue $\lambda$ satisfies
%---------------------------------------------------
\begin{equation}
e^{-\lambda L}=\frac{\alpha_1b_{22}+\alpha_2b_{11}\pm
\sqrt{(\alpha_1b_{22}+\alpha_2b_{11})^2-4\alpha_1\alpha_2(b_{11}b_{22}-b_{12}b_{21})}}
{2(b_{11}b_{22}-b_{12}b_{21})},
\label{criteria}\end{equation}
%---------------------------------------------------
which leads to the following theorem.
%---------------------------------------------------
\begin{theorem} With $\overline{B}$ determined by
\eqref{lin} and \eqref{bbar}. Let $\Lambda_1, \Lambda_2$ denote the roots of the quadratic
equation:
\begin{equation}\alpha_1\alpha_2-\Lambda(\alpha_1b_{22}+\alpha_2b_{11})+
\Lambda^2(b_{11}b_{22}-b_{12}b_{21})=0.\label{quadratic}\end{equation}
Assume that $|\Lambda_i|>1$ for $i=1, 2$, then the steady state $\overline{u}(x)$
is linearly stable for the system \eqref{Model:ode}.
\end{theorem}
%---------------------------------------------------
%---------------------------------------------------
%\begin{proof}
{\bf{Proof}}
To confirm the linear stability of the steady state $\overline{u}$,
it will suffice to show that all eigenvalues $\lambda$ to the eigenvalue
problem \eqref{spectral} are located in the left half of the complex plane, i.e.
$\Re(\lambda) <0$. Then from the equation $e^{-L\lambda}=\Lambda_i, i=1, 2$,
we have all the eigenvalues
$\lambda = -\frac{\ln |\Lambda_i|}{L}+ i(\arg \Lambda_i +2k\pi)$ for all $i=1,2$
and $k\in \mathbb Z$.
Therefore $|\Lambda|>1$ implies the linear stability of $\overline{u}$. The proof
is completed. $\hfill\square$
%\end{proof}
%---------------------------------------------------

In fact, the equation \eqref{criteria} can play the role of criteria for
the spectral stability of steady state $\overline{u}$ obtained in
Theorem~\ref{steadySol}. It is straightforward to notice that in the case of two
complex conjugate solutions to the above quadratic equation \eqref{criteria}, i.e.,
\[
(\alpha_1b_{22}+\alpha_2b_{11})^2-4\alpha_1\alpha_2(b_{11}b_{22}-b_{12}b_{21}))\leq 0.
\]
we have
\[
|e^{-\lambda L}|=\sqrt{\frac{\alpha_1\alpha_2}{(b_{11}b_{22}-b_{12}b_{21})}}.
\]
In the case of real roots, which happens to be the case for the numerical
experiment based on physiological parameters (with $L=15, c = 4.2$)
in Section~\ref{Sec:sim}, we have
%---------------------------------------------------
\[\left[\begin{array}{cc}
b_{11} & b_{12}\\
b_{21} & b_{22}\end{array}\right]
=
\left[\begin{array}{cc}
0.9978 & -0.0004\\
1.9048 & 0.8068\end{array}\right].
\]
%---------------------------------------------------
The quadratic equation \eqref{quadratic} becomes
\[2-2.8054\Lambda+0.8057\Lambda^2=0,\]
which has roots $\Lambda_1= 1.0003$ and $\Lambda_2 =2.4817$. As an immediate consequence, we know
that the steady solution $\overline{u}$ to the system \eqref{Model:ode} is linearly stable. Computations for other cases in Section~\ref{Sec:sim} reveal the stability as well.

%=======================================================================
%***********************************************************************
\subsection{Paritial differential equation $\varepsilon\neq0$ \label{subsect:epsilon}}
It is our belief that the second order diffusions have ignorable effects
on the distribution of glucose and insulin inside the pancreas. In this
section we justify this claim in the context of singular perturbation techniques.
It should be pointed out that our analysis on the boundary condition is novel
in the literature in the sense that the corresponding equations have mixed boundary
and initial conditions, the existing perturbation theory is not applicable to our
conditions under consideration.
In this section we discuss the steady state solutions of \eqref{Model_pancreas}
%---------------------------------------------------
\begin{equation}
\left\{\begin{aligned}
&-\varepsilon \frac{\partial^2 G}{\partial x^2} + \frac{\partial G}{\partial x}=G_{in}-a G I,\\
&-\varepsilon \frac{\partial^2 I}{\partial x^2} + \frac{\partial I}{\partial x}
=\frac{\sigma G^2}{b^2+G^2}-d_i I,\\
&G(L)=\alpha_1 G(0),\qquad I(L)=\alpha_2 I(0),\\
&G^\prime(L) = G^\prime(0),\qquad I^\prime(L)=I^\prime(0).
\end{aligned}\right.\label{steadywithdiffusion}\end{equation}
%---------------------------------------------------
Introducing that $U_\varepsilon = u = (G, I)^t \mbox{ and } V_\varepsilon = U^\prime_\varepsilon$,
we have the system
%---------------------------------------------------
\begin{equation}\label{steadyEquwithdiffusion}
-\varepsilon\frac{\partial^2 U_\varepsilon}{\partial x^2}+
\frac{\partial U_\varepsilon}{\partial x}
= B(U_\varepsilon -u^\star)+F(U_\varepsilon(x)).
\end{equation}
%---------------------------------------------------
Equivalently
%---------------------------------------------------
\[\left(\begin{array}{c}U_\varepsilon(x)\\ V_\varepsilon(x)\end{array}\right)^\prime
={\mathbb B}\left(\begin{array}{c}U_\varepsilon(x)-u^\star\\ V_\varepsilon(x)\end{array}\right)
-\left(\begin{array}{c} 0\\ \frac{1}{\varepsilon}F(U_\varepsilon)\end{array}\right),
\]
%---------------------------------------------------
with initial condition $U_\varepsilon(0) = u_0$ and $V_\varepsilon(0) = v_0$. In this section
we are going to search for appropriate $u_0, v_0$ so that $U_\varepsilon, V_\varepsilon$ satisfy the
the boundary conditions in \eqref{steadyEquwithdiffusion}.

We take the limit as $\varepsilon\rightarrow 0$, formally the model
\eqref{steadywithdiffusion} turns out to be
%---------------------------------------------------
\begin{equation}
\begin{cases}
\displaystyle\frac{\partial G}{\partial x}=G_{in}-a G I,\vspace{0.2cm}\\
\displaystyle\frac{\partial I}{\partial x}
=\frac{\sigma G^2}{b^2+G^2}-d_i I,\vspace{0.2cm}\\
G(L)=\alpha_1 G(0),\qquad I(L)=\alpha_2 I(0),
\end{cases}\label{steadywithoutdiffusion}\end{equation}
%---------------------------------------------------
where the first order derivative boundary conditions in \eqref{steadyEquwithdiffusion} are ignored.
The model with diffusion \eqref{steadywithdiffusion} should be regarded singular perturbation of \eqref{steadywithoutdiffusion} because the highest derivative in \eqref{steadywithdiffusion} vanishes in the limiting process as $\varepsilon\rightarrow 0$. As is well known in the singular perturbation theory, boundary layer or fast oscillation phenomenon might occur which prevent the convergence of \eqref{steadywithdiffusion} to the solution of \eqref{steadywithoutdiffusion} as $\varepsilon\rightarrow 0$. In this section, we claim that neither of them will appear in the above approximation. Toward this end, we first present the comparison between the two linearization problems with boundary value conditions:
%---------------------------------------------------
\begin{equation}
\begin{split}
&-\varepsilon\frac{\partial^2 u}{\partial x^2} + \frac{\partial u}{\partial x}= Bu,\\
\\
&u^\prime(L) = u^\prime(0),\qquad u(0) = u_0,
\end{split}\label{linearizationwithdiffusion}\end{equation}
%---------------------------------------------------
and
%---------------------------------------------------
\begin{equation}
\frac{\partial u}{\partial x}= Bu,\qquad u(0) = u_0.
\label{linearizationnodiffusion}\end{equation}
%---------------------------------------------------
We denote the solution of \eqref{linearizationwithdiffusion}
by $u=\Phi(x)u_0$, where $\Phi(x)$ will be defined in the lemma below,
while the solution of \eqref{linearizationnodiffusion} can be expressed
as $u = e^{Bx}u_0$.
%---------------------------------------------------
\begin{lemma}\label{comparison1}
For $0<<\varepsilon<1$, we have
\begin{equation}
\begin{split}
&\Phi(x)u_0-e^{Bx}u_0 = O(\varepsilon) u_0,\\
&\Phi^\prime(x)u_0-Be^{Bx}u_0 = O(1) u_0,
\end{split}\label{comp1}\end{equation}
where $O(\cdot)$ is uniform with respect to the size of $u_0$ and $\varepsilon$.
\end{lemma}
The proof of Lemma \ref{comparison1} can be found in Appendix \ref{Appen:Lem_Comparison1}. For the sake of convenience,  we introduce some a priori estimates on the solution to the non-homogeneous system
with singular parameter $\varepsilon>0$.
\begin{equation}
\begin{split}
&-\varepsilon\frac{\partial^2 u}{\partial x^2}
+\frac{\partial u}{\partial x}= Bu+g(x),\\
&u^\prime(L) = u^\prime(0),\qquad u(0) = 0.
\end{split}\label{diffusionwithg}\end{equation}
and
\begin{equation}
\frac{\partial u}{\partial x}= Bu+g(x),\qquad u(0) = 0.
\label{nodiffusionwithg}\end{equation}
We denote the solution of \eqref{diffusionwithg} by $u = [Sg](x)$ while the
solution of \eqref{nodiffusionwithg} can be expressed as
$u = \int^x_0e^{B(x-y)}g(y)dy$.

\begin{lemma}\label{comparison2}
For $0<< \varepsilon <1$,
\begin{equation}\begin{split}
[Sg](x) &-\int^x_0e^{B(x-y)}g(y)dy = O(\varepsilon)\int^x_0e^{B(x-y)}g(y)dy\\
&+O(1)\int^L_xe^{D_\varepsilon(x-y)}g(y) dy
+O(1)e^{D_0 x}\int^L_0e^{-D_\varepsilon y}g(y)dy\\
&+O(\varepsilon) e^{-D_\varepsilon(L-x)}\int_0^x e^{D_0(L-y)}g(y) dy,
\end{split}\label{comp2}\end{equation}
where $O(\cdot)$ is uniform with respect to the size of $u_0$ and $\varepsilon$.
\end{lemma}
The proof of the above lemma is presented in Appendix \ref{Appen:Lem_Comparison2}. Thanks to Duhammel principle, we can write the solution of system \eqref{steadywithdiffusion} as an integral equation
%---------------------------------------------------
\[\begin{split}
U_\varepsilon (x) & = ({\bf I}-D_{11}(x))u^\star
+D_{11}(x)u_0+D_{12}(x)v_0-\frac{1}{\varepsilon}\int_0^xD_{12}(x-y)F(U_\varepsilon(y))dy,\\
V_\varepsilon (x) & =
-D_{21}(x)u^\star+D_{21}(x)u_0+D_{22}(x) v_0
-\frac{1}{\varepsilon}\int_0^xD_{22}(x-y)F(U_\varepsilon(y))dy.
\end{split}\]
%---------------------------------------------------
By utilization of the boundary condition $V_\varepsilon(0)=V_\varepsilon(L)$,
we obtain an equation on $v_0$:
\[\begin{split}
v_0 &= -({\bf I}-D_{22}(L))^{-1}D_{21}(L)u^\star+({\bf I}-D_{22}(L))^{-1}D_{21}(L)u_0\\
&-\frac{1}{\varepsilon}(({\bf I}-D_{22}(L))^{-1}\int^L_0D_{22}(L-y)F(U_\varepsilon(y))dy.
\end{split}\]
Plugging $v_0$ into $U_\varepsilon(x)$leads to the integral equation on $U_\varepsilon(x)$:
\[\begin{split}
U_\varepsilon(x) & = ({\bf I} - \Phi(x))u^\star+\Phi(x) u_0\\
&-\frac{1}{\varepsilon}D_{12}(x) ({\bf I}-D_{22}(L))^{-1}
\int^L_0D_{22}(L-y)F(U_\varepsilon(y))dy\\
&-\frac{1}{\varepsilon}
\int^x_0D_{12}(x-y)F(U_\varepsilon(y))dy.\end{split}\]
Furthermore by the boundary condition $U_\varepsilon(L) =DU_\varepsilon(0)$,
an closed integral equation on $U_\varepsilon(x)$ shall be derived
\[\begin{split}
u_0 & = (D-E(L))^{-1}({\bf I} - E(L))u^\star\\
&-\frac{1}{\varepsilon}(D-E(L))^{-1}D_{12}(L) ({\bf I}-D_{22}(L))^{-1}
\int^L_0D_{22}(L-y)F(U_\varepsilon(y))dy\\
&-\frac{1}{\varepsilon}(D-E(L))^{-1}
\int^L_0D_{12}(L-y)F(U_\varepsilon(y))dy,\end{split}\]

\begin{equation}\begin{split}
U_\varepsilon(x)&= u^\star+\Phi(x)(D-E(L))^{-1}({\bf I}-D)u^\star\\
&-\frac{1}{\varepsilon}\Phi(x)(D-E(L))^{-1}D_{12}(L) ({\bf I}-D_{22}(L))^{-1}
\int^L_0D_{22}(L-y)F(U_\varepsilon(y))dy\\
&-\frac{1}{\varepsilon}D_{12}(x)({\bf I}-D_{22}(L))^{-1}
\int^L_0D_{22}(L-y)F(U_\varepsilon(y))dy\\
&-\frac{1}{\varepsilon}\Phi(x)(D-E(L))^{-1}
\int^L_0D_{12}(L-y)F(U_\varepsilon(y))dy\\
&-\frac{1}{\varepsilon}\int^x_0D_{12}(x-y)F(U_\varepsilon(y))dy\\
&=U_\varepsilon^0(x) + U_\varepsilon^1(x)+ U^2_\varepsilon(x).
\end{split}\label{finaliteration}\end{equation}
We make some comments on $U_\varepsilon^i, i=0, 1, 2$. First
\[
U_\varepsilon^0(x)=u^\star+\Phi(x)(D-E(L))^{-1}({\bf I}-D)u^\star=u^\star+\xi_\var^\star(x),
\]
which satisfies the following system:
\begin{equation}\label{iteration1}
-\varepsilon\frac{\partial^2 u}{\partial x^2}
+\frac{\partial u}{\partial x} = B(u-u^\star),
\qquad u^\prime(0) = u^\prime(L), \quad u(L) = D u(0).
\end{equation}
and
\[\begin{split}
U^2_\varepsilon(x) &= [SF(U_\varepsilon(\cdot)](x)\\
&=-\frac{1}{\varepsilon}D_{12}(x)({\bf I}-D_{22}(L))^{-1}
\int^L_0D_{22}(L-y)F(U_\varepsilon(y))dy\\
&-\frac{1}{\varepsilon}\int^x_0D_{12}(x-y)F(U_\varepsilon(y))dy,
\end{split}\]
which satisfies the system of differential equations:
\begin{equation}\label{iteration2}-\varepsilon\frac{\partial^2 u}{\partial x^2}
+\frac{\partial u}{\partial x} = Bu +F(U_\varepsilon(y)),
\quad u^\prime(0) = u^\prime(L), \quad u(0)=0.\end{equation}
Lastly,
\[\begin{split}
U_\varepsilon^1(x) &= \Phi(x)(D-E(L))^{-1}[SF(U_\varepsilon(\cdot))](L)\\
=&-\frac{1}{\varepsilon}\Phi(x)(D-E(L))^{-1}D_{12}(L) ({\bf I}-D_{22}(L))^{-1}
\int^L_0D_{22}(L-y)F(U_\varepsilon(y))dy\\
&-\frac{1}{\varepsilon}\Phi(x)(D-E(L))^{-1}
\int^L_0D_{12}(L-y)F(U_\varepsilon(y))dy,
\end{split}\]
which satisfies the following equations:
\begin{equation}\label{iteration3}
-\varepsilon\frac{\partial^2 u}{\partial x^2}
+\frac{\partial u}{\partial x} = Bu,
\quad u^\prime(0) = u^\prime(L), \quad u(0)
=(D-E(L))^{-1}[SF(U_\varepsilon(\cdot))](L).\end{equation}
It is easy to check that $U_\varepsilon^1+U_\varepsilon^2$ satisfies
the boundary condition
\[U_\varepsilon^1(L)+U_\varepsilon^2(L)
=D(U_\varepsilon^1(0)+U_\varepsilon^2(0)).\]
In this sense $U_\varepsilon^2$ is a modification on the boundary value of
$U_\varepsilon^1$ so that the boundary condition $u(L) = Du(0)$ is secured.
Similarly in the expression of $U_\varepsilon^2$ the first term is the modification
of the boundary value of the second term so that the boundary condition
$u^\prime(0)=u^\prime(L)$ is guaranteed.

Similar to the existence of steady state solution to \eqref{steadyEqua},
we define the nonlinear operator $T_\varepsilon$ from $C([0,L])^2$ to itself
by setting $T_\varepsilon U_\varepsilon$ as the right hand side in
the above equality. Then our task is to prove the existence of fixed point of nonlinear
operator $T_\varepsilon$ in $C([0,L])^2$. To this end we introduce the ball
$B_{r_\varepsilon}(u^\star)$, then prove that
$T_\varepsilon B_{r_\varepsilon}(u^\star)\subset B_{r_\varepsilon}(u^\star)$.
We are in a position to prove the main results in this section.
%---------------------------------------------------
\begin{theorem}\label{SSwithdiffusion}Assume that
\[
4\|\xi_\varepsilon^\star\|_\infty Lk(\|(D-E(L))^{-1}\|_2+1)\leq 1.
\]
Let
\[
r_\varepsilon = \frac{1-\sqrt{1-4\|\xi_\varepsilon^\star\|_\infty
Lk(\|(D-E(L))^{-1}\|_2+1)}}{2}.
\]
Then, for any $U\in B_{r_\varepsilon}(u^\star)$,
\[
\|T_\varepsilon U-u^\star\|_\infty\leq \|(D-E(L))^{-1}({\bf I}-Du^\star)\|_2
+Lk(\|(D-E(L))^{-1}\|_2+1)r^2_\varepsilon\leq r_\varepsilon,
\]
there exists a steady solution, denoted by $\overline{U}_\varepsilon(x)$, to the evolution
equation (\ref{Model_pancreas}), which satisfies
\[
\|\overline{U}_\varepsilon-u^\star\|_\infty\leq r_\varepsilon.
\]
Moreover we have the convergence $\overline{U}_\varepsilon\rightarrow \overline{u}$
in the $L^\infty$ norm as $\varepsilon\rightarrow 0$, where the $\overline{u}$ is the steady
state obtained in Theorem \ref{steadySol}.
\end{theorem}
%---------------------------------------------------
%---------------------------------------------------
%\begin{proof}
{\bf{Proof}}
The proof can be carried out along the same approach as the proof of Theorem \ref{steadySol} in observing of Lemma \ref{comparison1} and Lemma \ref{comparison2}. So we omit the details. $\Box$
%\end{proof}
%---------------------------------------------------
\section{Numerical analysis} \label{Sec:sim}
%---------------------------------------------------

Our main aim of this paper is to investigate and estimate the insulin distribution in pancreatic vein through mathematical modeling approach due to the lack of technology. However on the other hand, our estimation cannot be verified through {\em in vivo} experiments with current available technology. Thus we carefully choose physiologically reasonable parameter values identified in the literature and well known physiological facts and anatomic facts. Table \ref{Table:params} lists the ranges of these values and the references from which the ranges are determined.
%------------------------------------------
\begin{table} [!htbp] %[t]
\begin{center} \caption{\hskip 0.1in Ranges of parameter values used for numerical studies of Model (\ref{Model_pancreas}).} %\vspace{.3cm}
%\smallskip
\begin{tabular}{ c c l | l } %r r r | r r r }
\hline Parameters & Values & Units & References \\
\hline
\hline $c$ & 0.5--9 & cm/min & \citet[Fig. 5C]{nyman2010glucose} \\
\hline $b$ & 9 & mM & \citet{de2000mathematical, li2012range} \\
\hline $L$ & 15 & cm &  Estimated by anatomy \\
\hline $G_{in}$ & 0.06 & mM & \citet{li2006modeling, sturis1991computer} \\
\hline $a$ & $10^{-5}$ & mM/pM &  \citet{de2000mathematical, li2012range} \\
\hline $\sigma$ & 15 & pM/mM &  \citet{li2012range, wang2013dynamic} \\
\hline $d_i$ & 0.04 & min$^{-1}$ &  \cite{li2006modeling} \\
\hline $\alpha_1$ & 1 &   &  well known physiological fact \\
\hline $\alpha_2$ & 2 &   &  well known physiological fact \\
\hline
\end{tabular}
\label{Table:params}
\end{center}
\end{table}
%------------------------------------------

For the velocity of blood flow $c$, as discussed at the beginning of Section \ref{Sec:sssol}, the velocity of blood flow is assumed to be 1 cm/min for simplicity in this analytical study. In numerical studies, we determine the velocities according to \cite{nyman2010glucose} as follows. For the hypoglycemic case, $c$ is in the range from 50 to 500 $\mu$m/s with mean at 100 $\mu$m/s, which can be converted to the range from 0.5 cm/min to 3 cm/min. For the case of hyperglycemia, the range is from 100 to 1500 $\mu$m/s with mean at 700 $\mu$m/s, or, equivalently, 3 cm/min to 9 cm/min with mean at 4.2 cm/min.

Theorem~\ref{steadySol} and  Theorem~\ref{SSwithdiffusion} provide not only the theoretical existence for the steady solutions with or without small diffusion, but also the numerical scheme
for tracking the steady solutions through iterations. We use Theorem~\ref{SSwithdiffusion}
as an example. Initially we set up $u^0 = U^0_\varepsilon$, which is corresponding to solving
differential equation \eqref{iteration1}, and $u^{n+1} = T_\varepsilon u^{n}$, which equivalently
amounts to solving the differential equations \eqref{iteration1} and \eqref{iteration2}.

We first investigate how different boundary distributions of insulin affect the steady state distribution of insulin. We consider five simplified and possibly true boundary distributions of insulin: homogeneous distribution thorough out the vein;  increasingly linear distributed;  decreasingly linear distributed; concave upward quadratic distribution; and concave downward quadratic distribution. Fig \ref{Fig_NoDiff_Homo}--\ref{Fig_NoDiff_RevQuadLinear} show the steady state solutions for each case. A widely believed hypothesis is that insulin concentration inside pancreas is much higher than concentrations of insulin in periphery and portal vein. However, according to our numerical studies, this hypothesis is only possibly true when the boundary distribution of insulin secreted into the tail area of the pancreatic vein is higher, which may implies that higher $\beta$-cell mass in the tail area results in higher insulin concentration in most part of pancreatic vein. On the other hand, if the above hypothesis is true, we can inversely induce that more $\beta$-cells reside in the tail area inside pancreas than the head area. See Fig. \ref{Fig_NoDiff_DecLinear} and Fig \ref{Fig_NoDiff_RevQuadLinear}.

%------------------------------------------
\begin{figure} [!htbp]
	\centerline{\resizebox{2.5in}{1.5in}{\includegraphics{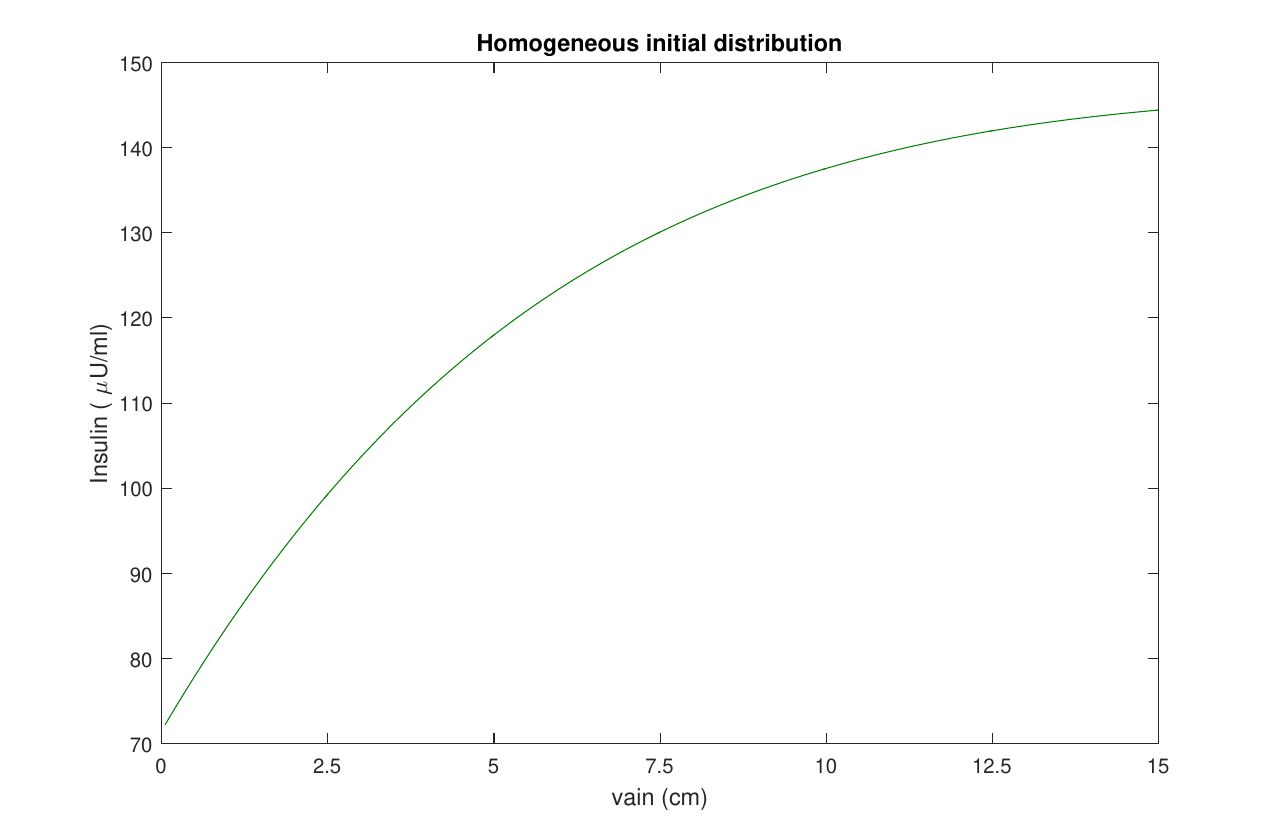}}
		\resizebox{2.5in}{1.5in}{\includegraphics{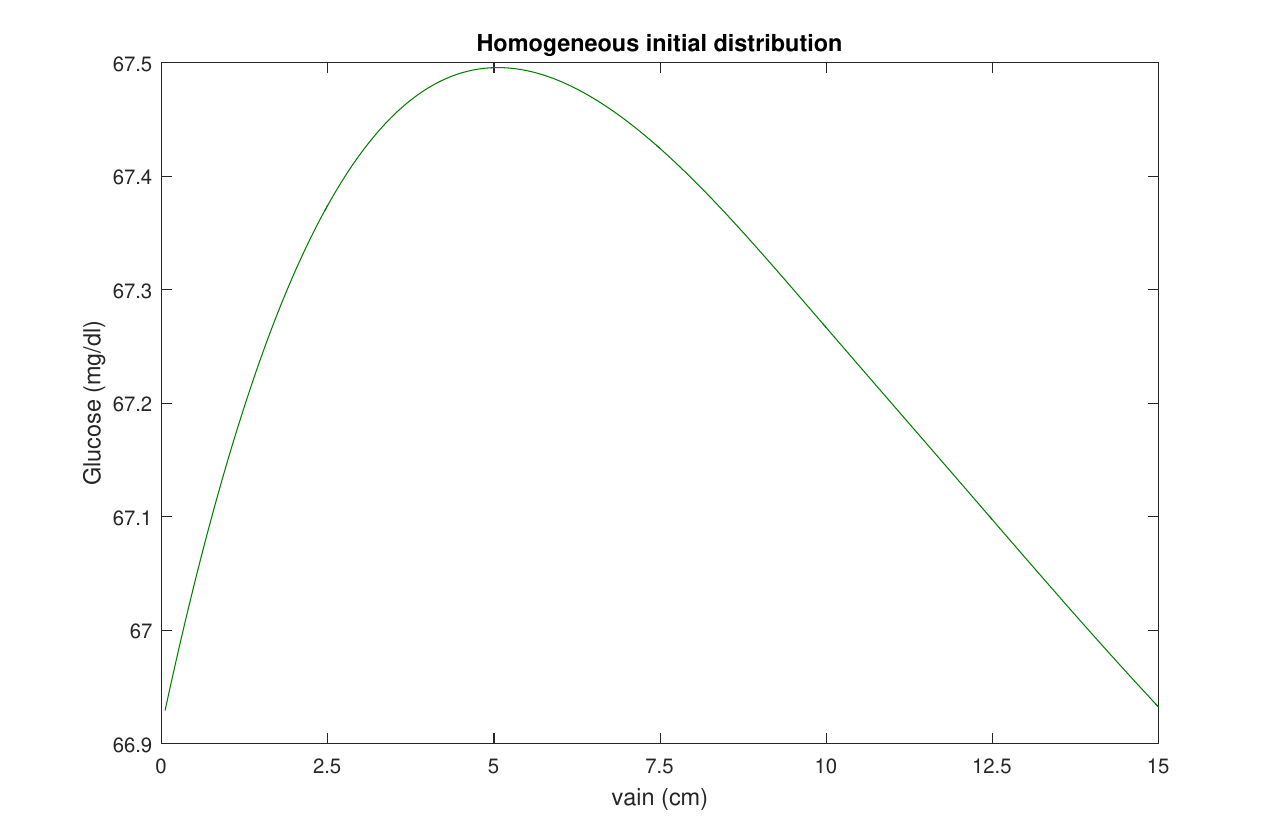}}}
	\caption{\footnotesize \hskip 0.2in Distribution of insulin and glucose in pancreatic vein for homogeneous input.}
	\label{Fig_NoDiff_Homo}
\end{figure}
%------------------------------------------

%------------------------------------------
\begin{figure} [!htbp]
	\centerline{\resizebox{2.5in}{1.5in}{\includegraphics{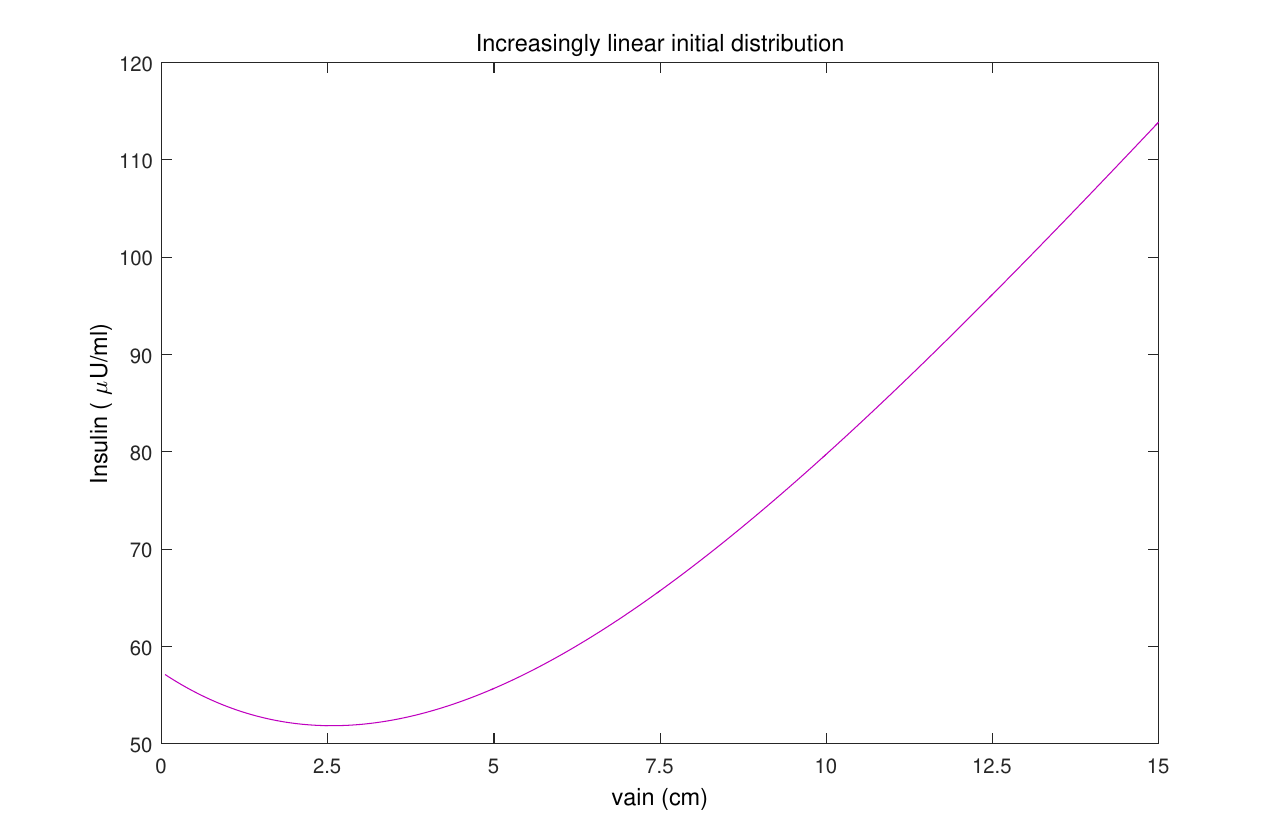}}
		\resizebox{2.5in}{1.5in}{\includegraphics{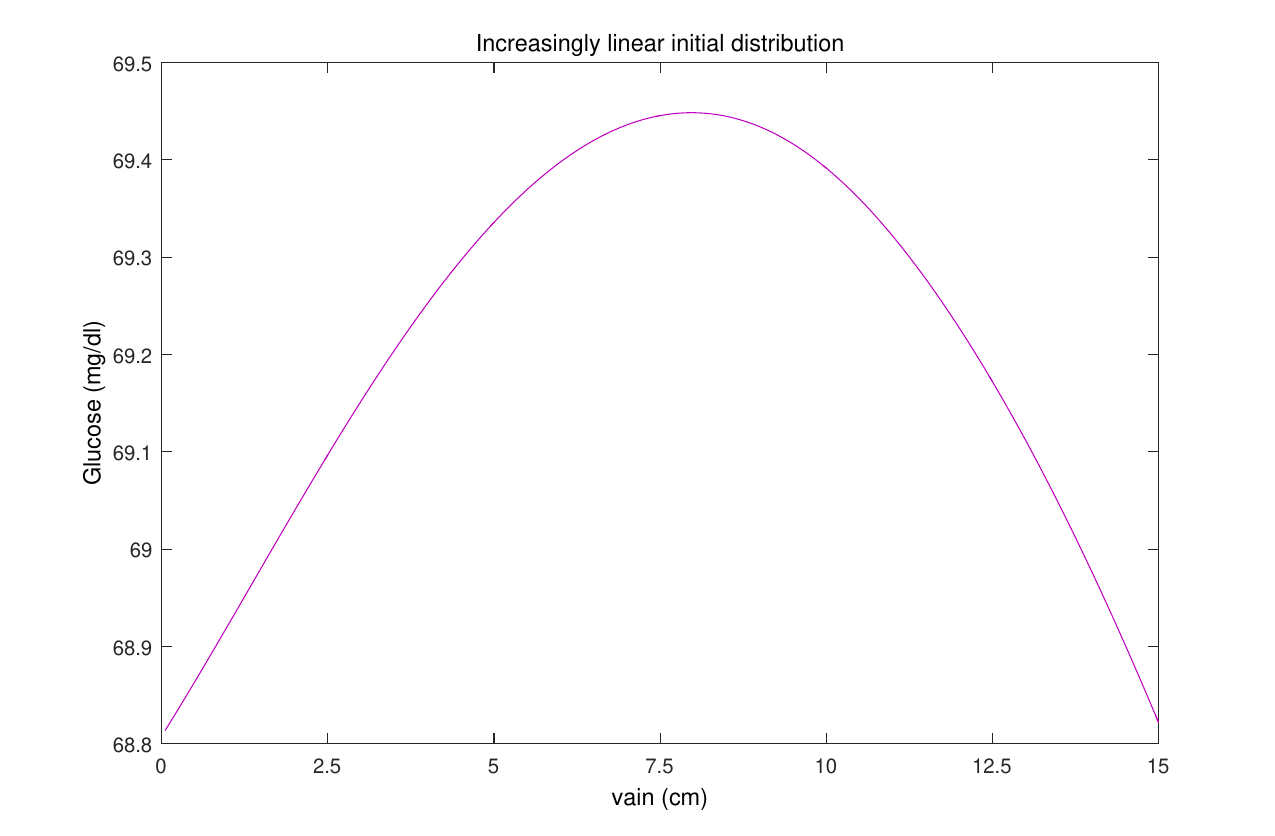}}}
	\caption{\footnotesize \hskip 0.2in Distribution of insulin and glucose in pancreatic vein for increasing input.}
	\label{Fig_NoDiff_Linear}
\end{figure}
%------------------------------------------

%------------------------------------------
\begin{figure} [!htbp]
	\centerline{\resizebox{2.5in}{1.5in}{\includegraphics{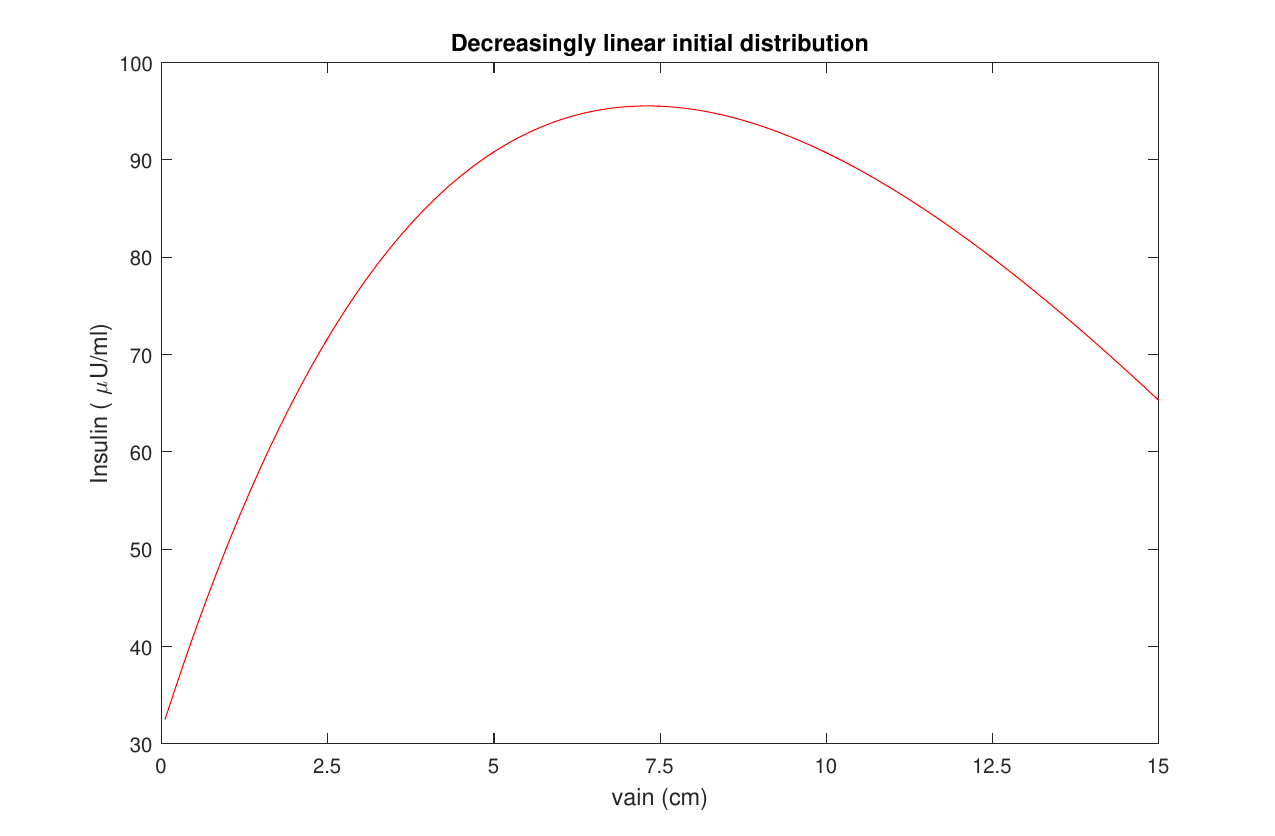}}
		\resizebox{2.5in}{1.5in}{\includegraphics{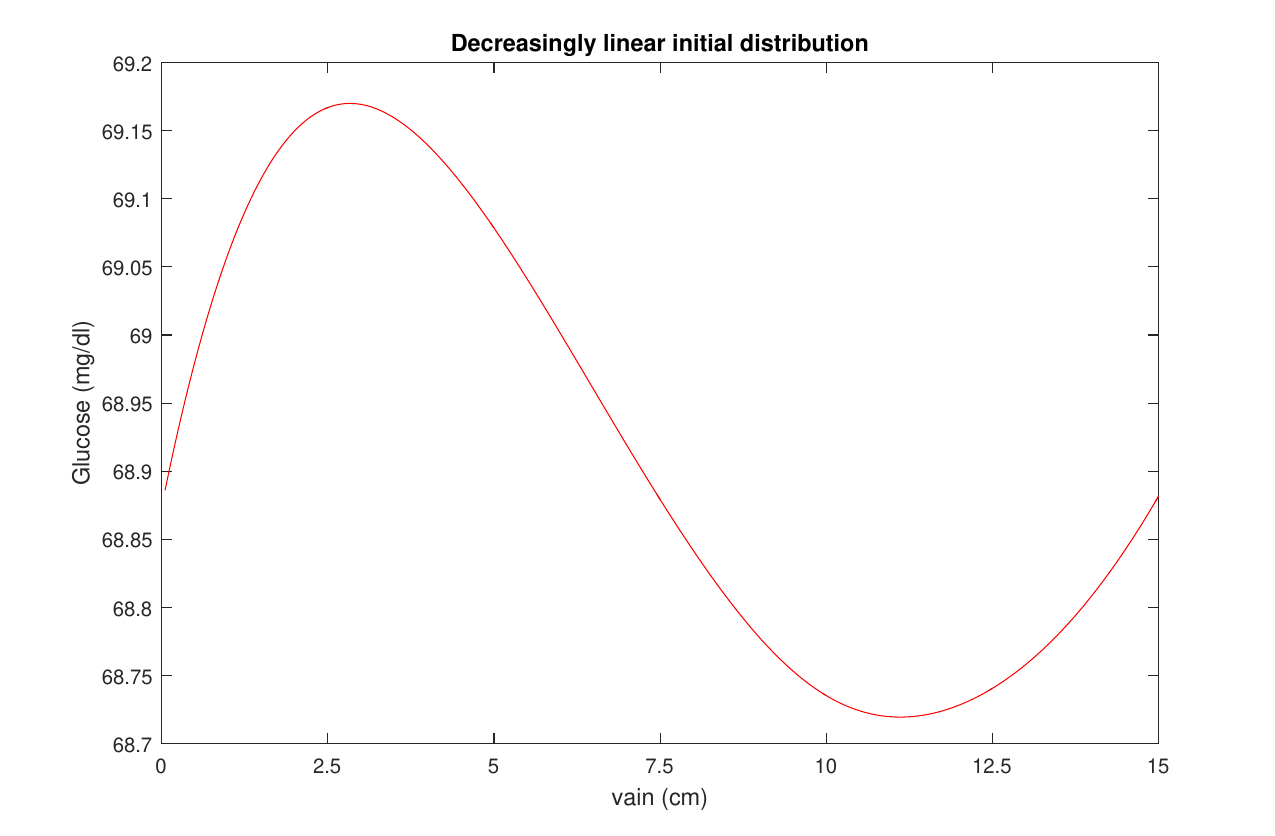}}}
	\caption{\footnotesize \hskip 0.2in Distribution of insulin and glucose in pancreatic vein for decreasing linear input.}
	\label{Fig_NoDiff_DecLinear}
\end{figure}
%------------------------------------------

%------------------------------------------
\begin{figure} [!htbp]
	\centerline{\resizebox{2.5in}{1.5in}{\includegraphics{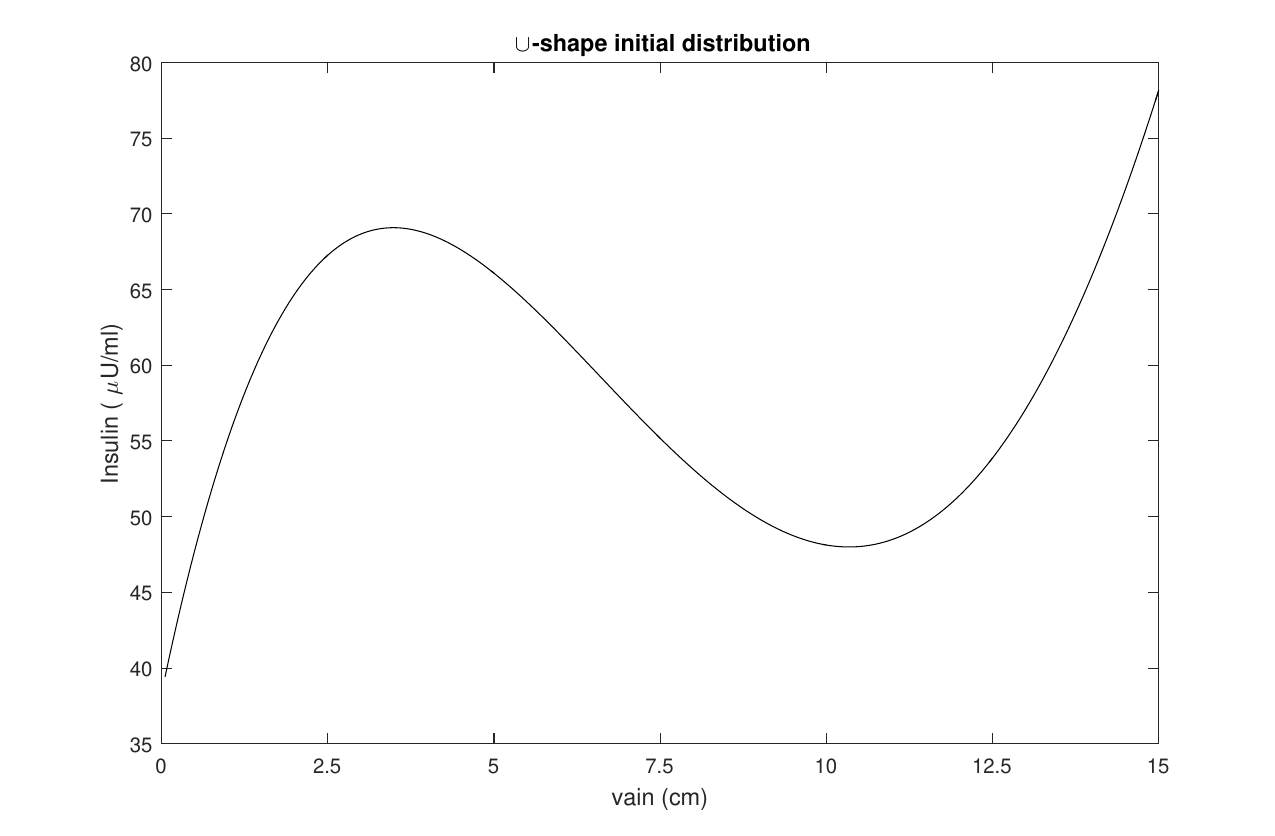}}
		\resizebox{2.5in}{1.5in}{\includegraphics{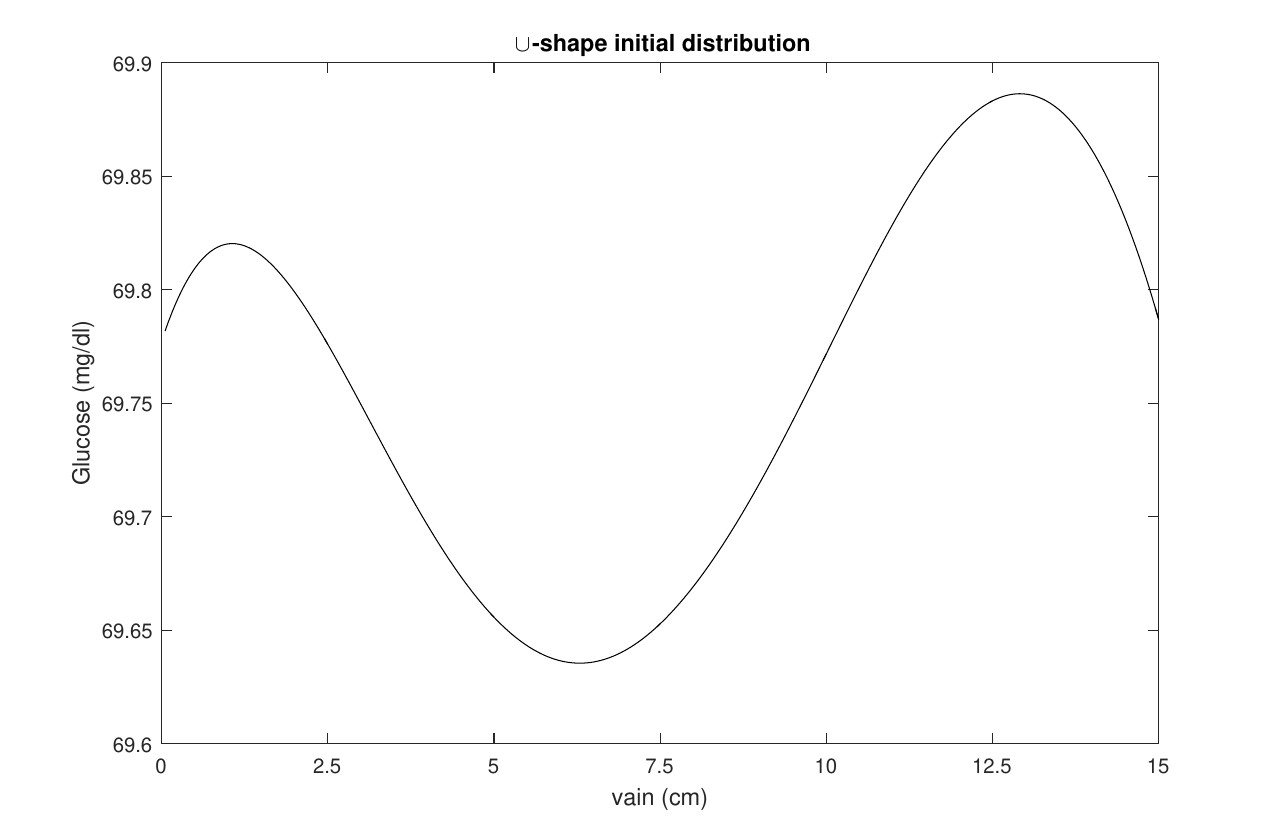}}}
	\caption{\footnotesize \hskip 0.2in Distribution of insulin and glucose in pancreatic vein for a quadratic input.}
	\label{Fig_NoDiff_QuadLinear}
\end{figure}
%------------------------------------------

%------------------------------------------
\begin{figure} [!htbp]
	\centerline{\resizebox{2.5in}{1.5in}{\includegraphics{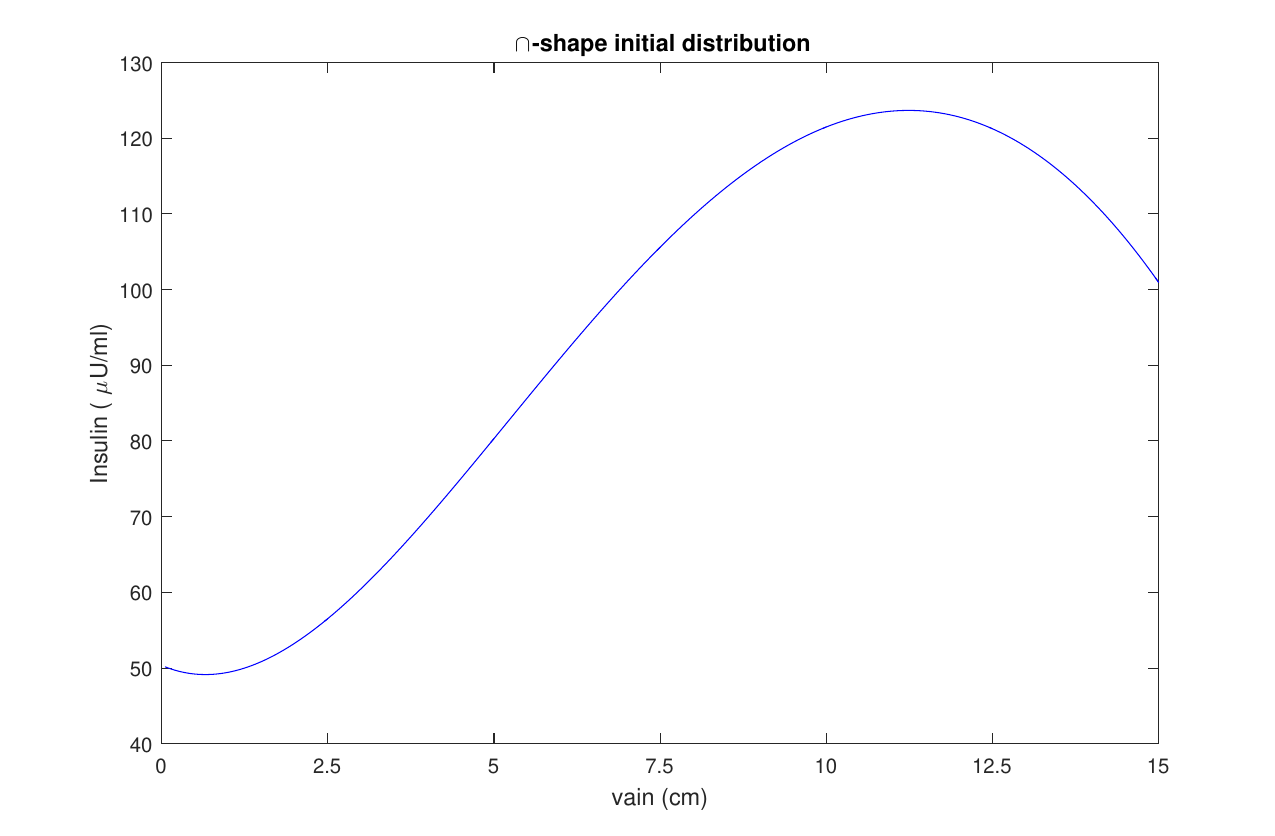}}
		\resizebox{2.5in}{1.5in}{\includegraphics{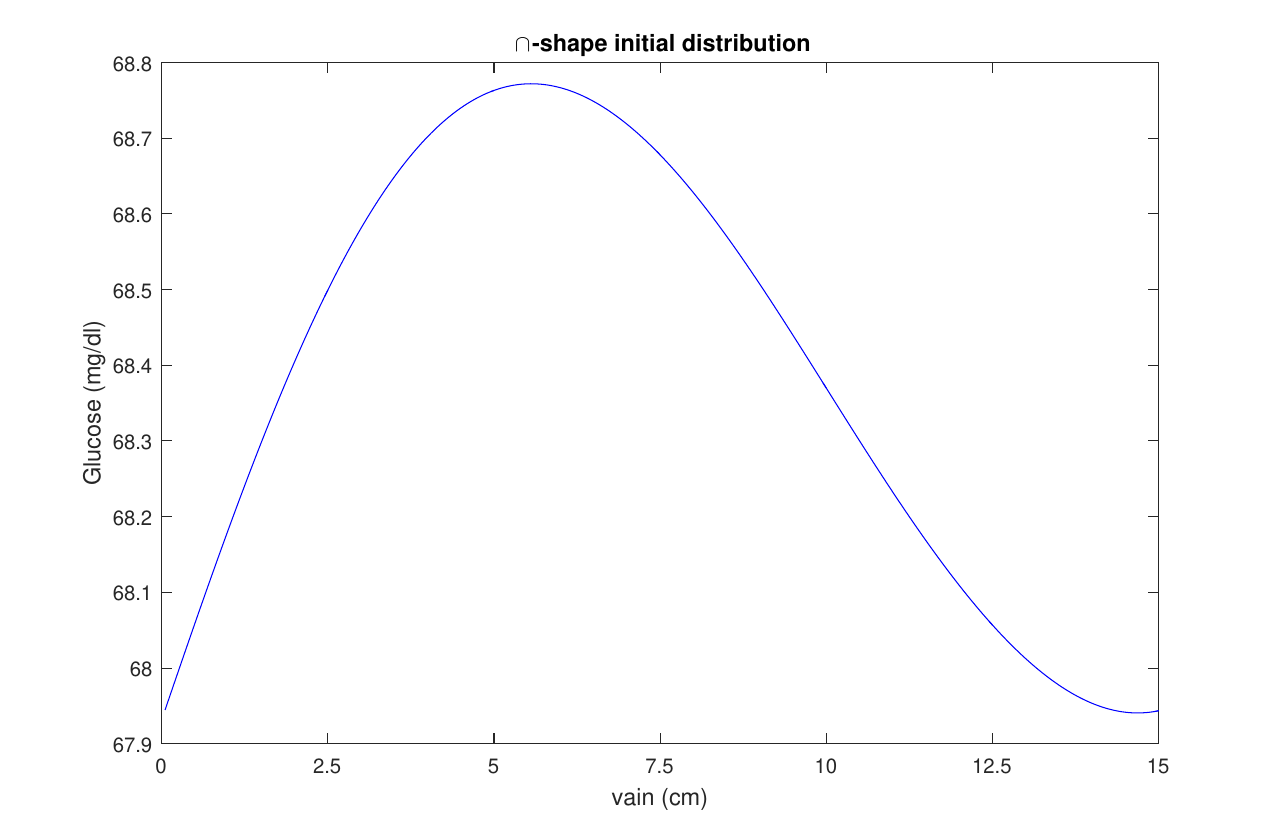}}}
	\caption{\footnotesize \hskip 0.2in Distribution of insulin and glucose in pancreatic vein for a reversed quadratic input.}
	\label{Fig_NoDiff_RevQuadLinear}
\end{figure}
%------------------------------------------

%\newpage

%{\color{red}
We will consider different velocity $c$ and some boundary distributions affect steady state distribution. For example,
We simulated various cases for the velocity of blood flow in pancreas according to the measurements in \cite{nyman2010glucose}. For hyperglycemic case, the velocity is as fast up to $c=9$ cm/min with mean value at $c=4.2$ cm/min. For hypoglycemic case, the velocity can be as slow as $c=0.5$ cm/min while the mean value is at $c=3$ cm/min.
%In Fig. \ref{Fig_L15_c0p5}, 
We noticed, when $c=0.5$ cm/min, glucose concentration is at about 3.3 mM while insulin concentration increases along the pancreas from tail to head. The increase is quicker in the tail and slower in the head area in the range about 112 pmol and 220 pmol.

When the velocity is at the average for hypoglycemic case,% (see Fig. \ref{Fig_L15_c3}),
$c=0.3$ cm/min, glucose level is in the range of 6 mM. Insulin concentration increases in similar manner but near linearly from about 65 pmol to 128 pmol.

In contrast, in Fig. \ref{Fig_L15_c4p2}, when $c=4.2$ cm/min, the average velocity in the hyperglycemic situation, the glucose level at 7 mM while insulin concentration increases nearly linearly from 56 pmol to 112 pmol. When the velocity is $c=9$ cm/min for extrem hyperglycimea, the glucose contration is at about 9.8 mM with nearly linearly increasing insulin concentration from 40 pmol to 80 pmol.% (see Fig. \ref{Fig_L15_c9}).

We can also calculate the diffusion rate $\varepsilon$ from \cite{nyman2010glucose}. However, since we have proved that a small diffusion rate does not disturb the dominance of convection, we omit numerical simulations.
}

%---------------------------------------------------
%================================================================

%================================================================
%---------------------------------------------------
\section{Discussions} \label{Sec:disc}
%---------------------------------------------------

In the present work we report on our mathematical modelling studies of the spatial distribution of insulin concentrations in pancreas using a novel PDE model.

We showed analytically that a small diffusion contribution $\varepsilon>0$ does not contribute to the spatial distribution of glucose and insulin along the pancreas. So our focus is on the case $\varepsilon = 0$. In this case, analytically the steady state solution exists and it is stable. Numerically, with physiologically reasonable parameters obtained from the literature, we demonstrated profiles of the distributions of both insulin and glucose. The simulations reveal that when glucose level is lower, which implies that the velocity of blood flow is slower, insulin concentration increases along pancreas from quicker to slower. When glucose level is higher, the velocity is faster and insulin concentration level increases nearly linearly. This is in agreement with the findings in \cite{nyman2010glucose} for blood flow velocities vs. the glucose concentration levels.

Due to the current lack of sufficient technology, such distributions cannot be tested by experiments. However, the profiles demonstrates reasonable and interesting distributions in agreement with the experiments performed by \cite{nyman2010glucose} and other basic physiological facts, for example, insulin concentration level is doubled at the head to the tail. It reveals that the length of pancreas and the velocity of blood flux likely play important roles for the distribution of insulin in pancreas.

%\begin{subappendices}
%	\subsection{How I became inspired}
%	What are here?
%\end{subappendices}

%---------------------------------------------------
%******************************************************************
%\appendixpage\appendix\section{Proof of main results\label{sect:proof}}
%-------------------------------------------------------------------
%----------------------------------------------------------------
%================================================================
\begin{appendices}
%----------------------------------------------------------------
\section{Proofs of results }
	
%----------------------------------------------------------------
\subsection{\bf Proof of Lemma \ref{lem:equiexistence}}\label{Appen:Lem_EquiExistence}
The uniqueness of the steady state can be seen from the
monotonicity of the right hand side of \eqref{Model:ode} with respect to
variables $G$ and $I$. The bounds for $I^\star$ is straightforward.
Bounds for $G^\star$ can be derived from the null-cline of the system
\eqref{Model:ode}.

To study the stability of $u^\star= (G^\star, I^\star)$  of system
\eqref{Model:ode},  we need to compute the eigenvalue of $B$,
%---------------------------------------------------
\[\mbox{det}|\lambda  -B|=0,\]
%---------------------------------------------------
which leads to the equation on $\lambda$,
%---------------------------------------------------
\[(\lambda+aI^\star)(\lambda+d_i)+\frac{2ad_i^2b^2I^{\star 2}}{\sigma G^{\star 2}}=0.\]
%---------------------------------------------------
then
%---------------------------------------------------
\[\lambda_{1,2}=\frac{-(aI^\star+d_i)\pm\sqrt{(aI^\star-d_i)^2
-\frac{2ad_i^2b^2I^{\star 2}}{\sigma G^{\star 2}}}}{2}.\]
%---------------------------------------------------
Since apparently $\Re(\lambda_{1,2}) <0$, so the equilibrium $u^\star=(G^\star,I^\star)$ is stable.
An routine phase portrait analysis can show that the equilibrium $u^\star$ is globally
stable. $\hfill\square$
%---------------------------------------------------

%{\bf Proof of Lemma \ref{lem:existstate}}
\subsection{Proof of Lemma \ref{lem:existstate}}\label{Appen:Lem_ExistState}
If $\alpha_j>1(j=1,2)$ and $\sigma(x)>\underline\sigma>0$, it follows from \eqref{eqn:I} that the nonlinear operator $I[G](x)$ is a positive operator for all $x\in[0,L],$
So does the operator $G.$

Substituting $I[G]$ into \eqref{eqn:G} yields
\begin{align}\label{eqn:GG}
G(x)=G_0e^{-\int_0^xI[G](s)ds}+G_{in}\int_0^xe^{-\int_\tau^xI[G](s)ds}d\tau:=\mathcal T[G](x).
\end{align}
Now, we are going to prove the compactness of the operator $\mathcal T[G]$. To do that end, we firstly prove the boundedness of $I[G]$, for each $x\in[0,L]$, we have that
\[|I[G](x)|\le I_0+\frac{\bar\sigma}{d}:=M_I,\forall~x\in[0,L].\]
For $x,\hat x\in[0,L],$ without loss of generality we assume that $\hat x<x.$ Note that
\begin{align}
|\mathcal T[G](x)-\mathcal T[G](\hat x)|\le& G_0\left|e^{-\int_0^xI(\tau)d\tau}-e^{-\int_0^{\hat x}I(\tau)d\tau}\right|\notag
\\
&+G_{in}\left|\int_0^x e^{-a\int_\tau^xI(s)ds}d\tau-\int_0^{\hat x}e^{-a\int_\tau^{\hat x}I(s)ds}d\tau\right|\notag
\\
\le&M_I|x-\hat x|+G_{in}\int_{\hat x}^xe^{-a\int_\tau^xI(s)ds}d\tau\notag
\\
&+G_{in}\int_0^{\hat x}\left|e^{-\int_\tau^xI[G](s)ds}-e^{-\int_\tau^{\hat x}I[G](s)ds}\right|d\tau\notag
\\
\le&(M_I+G_{in}+M_IG_{in}L)|x-\hat x|
\end{align}
here we have used the fact $|e^{-a}-e^{-b}|\le|a-b|.$ Therefore, the operator $\mathcal T$ is compact. From \protect[Corollary 10.1,\cite{I_2017}], we have that there exists at least one endemic steady state.
$\hfill\square$
%-------------------------------------------------------------------------------------------

%--------------------------------------------------------------------------------------------
%\noindent{\bf Proof of Lemma \ref{lem:steadyuniq}}
\subsection{Proof of Lemma \ref{lem:steadyuniq}}\label{Appen:Lem_SteadyUniq}
 By way of contradiction, we assume that \eqref{steadyEqua} has two positive solution $G_1$ and $G_2$ and $G_1\neq G_2$. By equation \eqref{eqn:I}, we note that 
\begin{align}\label{eqn:IG}
|I[G_1](x)-I[G_2](x)|\le& \int_0^x\left|\frac{\sigma(s) G_1^2(s)}{b^2+G_1^2(s)}-\frac{\sigma G_2^2(s)}{b^2+G_2^2(s)}\right|\notag
\\
=&\int_0^x\left|\frac{\sigma(s) b^2(G_1^2(s)-G_2^2(s))}{(b^2+G_1^2)(b^2+G_2^2)}ds\right|\notag
\\
\le&\bar\sigma b^2L.
\end{align}
Moreover, from \eqref{eqn:G}, we have that for each $x\in[0,L]$
\begin{align}\label{eqn:IG}
|G_1(x)-G_2(x)|\le& G_0\left|e^{-\int_0^xI[G_1](s)ds}-e^{-\int_0^xI[G_2](s)ds}\right|\notag
\\
&+G_{in}\int_0^x\left|e^{-\int_\tau^xI[G_1](s)ds}-e^{-\int_\tau^xI[G_2](s)ds}\right|d\tau\notag
\\
\le&G_0\bar\sigma b^2L^2+G_{in}\bar \sigma b^2 L^3
\end{align}
Picking up $L\le\max\left\{\frac{\epsilon}{2\sqrt{G_0\bar\sigma b^2L^2}},\frac{\epsilon}{2\sqrt[3]{G_{in}\bar\sigma b^2}}\right\},$ where $\epsilon$ is an sufficiently small value, we conclude that
$|G_1(x)-G_2(x)|\le\epsilon.$ By the arbitrary choice of $\epsilon$, we have that $G_1=G_2$ for each $x\in[0,L].$ 
$\hfill\square$

%---------------------------------------------------------------------------------------
%\noindent{\bf Proof of Lemma \ref{comparison1}} 
\subsection{Proof of Lemma \ref{comparison1}} \label{Appen:Lem_Comparison1}
Firstly, we recall that there are two negative eigenvalues
of $B$, denoted by $\Lambda_i, i=1, 2$, the corresponding eigenvectors are $P_1, P_2$.
We let $P = [P_1\quad P_2]$. So with matrix $P$, we have
\[B=P\mbox{diag} (\Lambda_1, \Lambda_2)P^{-1},\]
and $e^{Bx}u_0$ can be written as
\[e^{Bx} = P\mbox{diag}(e^{\Lambda_1 x}, e^{\Lambda_2 x})P^{-1}.\]
Back to model \eqref{linearizationwithdiffusion}, introduce $u^\prime = v$,
then the model turns out to be the augmented system:
\[u^\prime = v, \qquad v^\prime =-\frac{1}{\varepsilon}Bu+\frac{1}{\varepsilon} v,\]
or equivalently
\[\left(\begin{array}{c} u\\ v\end{array}\right)^\prime
= {\mathbb B}\left(\begin{array}{c} u\\ v\end{array}\right),
\]
where
\[{\mathbb B}=\left(\begin{array}{cc} 0 & {\bf I}\\
										\\
                    -\frac{1}{\varepsilon}B & \frac{1}{\varepsilon}{\bf I}
                    \end{array}\right).\]
The characteristic equation of the above system is
\[\mbox{det}|\varepsilon\lambda^2-\lambda+B| = 0,\]
which leads to
\[\varepsilon \lambda^2-\lambda+\Lambda_i =0,\qquad i=1, 2.\]
We denote the eigenvalues of order $O(1)$ by
\[\lambda_i^0 = \Lambda_i+O(\varepsilon),\]
corresponding eigenvectors are $(P_i,\lambda_i^0P)^t$ and
eigenvalues of order $O(1/\varepsilon)$,
\[\lambda_i^\varepsilon=\frac{1}{\varepsilon}-\Lambda_i+O(\varepsilon),\]
corresponding eigenvectors are $(P_i,\lambda_i^0P_i)^t$. We introduce following diagonal
matrices, $D_0 = \mbox{diag}(\lambda_1^0, \lambda_2^0), \mbox{ and } D_\varepsilon = \mbox{diag}(\lambda_1^\varepsilon,\lambda_2^\varepsilon)$. Let
\[{\mathbb P}
=\left(\begin{array}{cccc}
P_1 & P_2 & P_1 & P_2\\
\lambda_1^0 P_1
&\lambda_2^0 P_2
&\lambda_1^\varepsilon P_1
&\lambda_2^\varepsilon P_2
\end{array}\right) =
\left(\begin{array}{cc} P & 0\\ 0 & P\end{array}\right)
\cdot
\left(\begin{array}{cc} {\bf I} & {\bf I}\\ D_0 & D_\varepsilon
\end{array}\right).
\]
A straightforward computation leads to
\[{\mathbb P}^{-1}=
\left(\begin{array}{cc} \frac{D_\varepsilon}{D_\varepsilon -D_0}
& -\frac{{\bf I}}{D_\varepsilon -D_0}\\
\\
-\frac{D_0}{D_\varepsilon -D_0}
& \frac{{\bf I}}{D_\varepsilon -D_0}\end{array}\right)
\cdot
\left(\begin{array}{cc} P^{-1} & 0\\
0 & P^{-1}
\end{array}\right).
\]
Then the solution can be written
\[\begin{split}
&e^{{\mathbb B}x}\left(\begin{array}{c}u_0\\v_0\end{array}\right)
={\mathbb P}
\left(\begin{array}{cc} e^{D_0 x} & 0\\ 0 & e^{D_\varepsilon x}\end{array}\right)
{\mathbb P}^{-1}
\left(\begin{array}{c}u_0\\v_0\end{array}\right)\\
&=\left(\begin{array}{cc} P & 0\\ 0 & P\end{array}\right)
\left(\begin{array}{cc} \frac{D_\varepsilon e^{D_0x}-D_0e^{D_\varepsilon x}}{D_\varepsilon -D_0}
& \frac{e^{D_\varepsilon x}-e^{D_0 x}}{D_\varepsilon -D_0}\\
\\
\frac{D_\varepsilon D_0(e^{D_0x}-e^{D_\varepsilon x})}{D_\varepsilon -D_0}
& \frac{D_\varepsilon e^{D_\varepsilon x}-D_0e^{D_0 x}}{D_\varepsilon -D_0}\end{array}\right)
\left(\begin{array}{cc} P^{-1} & 0\\ 0 & P^{-1}\end{array}\right)
\left(\begin{array}{c}u_0\\v_0\end{array}\right).
\end{split}\]
From the computation we introduce matrix blocks of $e^{{\mathbb B} x}$:
\[e^{{\mathbb B}x}
=\left(\begin{array}{cc}
D_{11}(x) & D_{12}(x)\\
D_{21}(x) & D_{22}(x)\end{array}\right).\]
The definitions of $D_{ij}(x)$ are self-evident, and the following identities hold
\[D_{21}(x) = D^\prime_{11}(x),\qquad D_{22}(x) = D^\prime_{12}(x).\]
Let
\[\left(\begin{array}{c} u_\varepsilon(x)\\ v_\varepsilon(x)\end{array}\right)
= e^{{\mathbb B}x}
\left(\begin{array}{c} u_\varepsilon(0)\\ v_\varepsilon(0)\end{array}\right).\]
Then by boundary condition $v_\varepsilon(0) = v_\varepsilon(L)$, and
$v_\varepsilon(0) = ({\bf I}-D_{22}(L))^{-1}D_{21}(L)u_\varepsilon(0)$, we have
\[u_\varepsilon(x) = [D_{11}(x)+D_{12}(x)({\bf I}-D_{22}(L))^{-1}D_{21}(L)]
u_\varepsilon(0) = \Phi(x) u_\varepsilon(0).\]
and
\[v_\varepsilon(x) = u^\prime_\varepsilon(x).\]
By the boundary constraint that $u^\prime(0) = u^\prime(L)$, we have the equation on $v_0$:
\[
v_0=\frac{D_\varepsilon D_0(e^{D_0 L}-e^{D_\varepsilon L})}
{D_\varepsilon({\bf I} - e^{D_\varepsilon L})-D_0({\bf I}-e^{D_0L})} u_0
=({\bf I}+O(e^{-D_\varepsilon L}))D_0 u_0.\]
Plugging it into the solution, we obtain
\[u_\varepsilon (x) = (e^{Bx}+O(\varepsilon))u_0,\]
and
\[v_\varepsilon (x) = (Be^{Bx}+O(1))u_0,\]
which together imply the comparison. $\hfill\square$

%---------------------------------------------------------------------------------------
%\noindent{\bf Proof of Lemma \ref{comparison2}}
\subsection{Proof of Lemma \ref{comparison2}}\label{Appen:Lem_Comparison2}
Without loss of generality we can assume that $P={\bf I}$.
Otherwise we can consider the linear transformation $P^{-1}u, P^{-1}v$.
First we introduce the equivalent augmented system:
\[\left(\begin{array}{c} U_\varepsilon(x)\\ V_\varepsilon(x)\end{array}\right)^\prime
={\mathbb B}\left(\begin{array}{c} U_\varepsilon(x)\\ V_\varepsilon(x)\end{array}\right)
-\left(\begin{array}{c} 0 \\ \frac{1}{\varepsilon}g(x)\end{array}\right).\]
We have
\[
V_\varepsilon(x)  = D_{22}(x)V_\varepsilon(0)
-\frac{1}{\varepsilon}\int_0^x D_{22}(x-y)g(y)dy.\]
and the boundary condition $V_\varepsilon(0) = V_\varepsilon(L)$ yields
\[V_\varepsilon(0)=-\frac{1}{\varepsilon}({\bf I} - D_{22}(L))^{-1}
\int_0^L D_{22}(L-y)g(y)dy.\]
Then $U_\varepsilon$ can be solved as
\[\begin{split}
U_\varepsilon(x) & = -\frac{1}{\varepsilon}D_{12}(x)({\bf I} - D_{22}(L))^{-1}
\int_0^L D_{22}(L-y)g(y)dy\\
&-\frac{1}{\varepsilon}\int^x_0 D_{12}(x-y) g(y)dy
=-I_1+I_2-I_3+I_4,
\end{split}\]
where
\[\begin{split}
I_1&= \frac{1}{\varepsilon}D_{12}(x)({\bf I} - D_{22}(L))^{-1}
\int_0^L \frac{D_\varepsilon e^{D_\varepsilon (L-y)}}{D_\varepsilon -D_0}g(y)dy,\\
I_2&= \frac{1}{\varepsilon}D_{12}(x)({\bf I} - D_{22}(L))^{-1}
\int_0^L \frac{D_0 e^{D_0 (L-y)}}{D_\varepsilon -D_0}g(y)dy,\\
I_3&=\frac{1}{\varepsilon}\int_0^x\frac{ e^{D_\varepsilon (x-y)}}{D_\varepsilon -D_0}g(y)dy,\\
I_4&=\frac{1}{\varepsilon}\int_0^x\frac{ e^{D_0 (x-y)}}{D_\varepsilon -D_0}g(y)dy.
\end{split}\]
We observe that
\[({\bf I}-D_{22}(L))^{-1} = -(D_\varepsilon -D_0)
D_\varepsilon^{-1}e^{-D_\varepsilon L}({\bf I}+O(e^{-D_\varepsilon L})).\]
We focus on the singular term which contains $e^{D_\var x}$.
Therefore
\[\begin{split}
I_1&=\frac{1}{\varepsilon}\frac{e^{D_\varepsilon x}-e^{D_0 x}}{D_\varepsilon -D_0}
({\bf I}+O(e^{-D_\varepsilon L}))(D_\varepsilon -D_0)D_\varepsilon^{-1}
e^{-D_\varepsilon L}\int_0^L\frac{D_\varepsilon e^{D_\varepsilon (L-y)}}{D_\varepsilon -D_0}g(y) dy\\
&=\frac{1}{\varepsilon}({\bf I}+O(e^{-D_\varepsilon L}))
({\bf I}-e^{-(D_\varepsilon -D_0) x})
\int_0^L\frac{ e^{D_\varepsilon (x-y)}}{D_\varepsilon -D_0}g(y) dy.
\end{split}\]
Then we have
\[\begin{split}I_1+I_3 &= -\frac{1}{\varepsilon}
\int_x^L\frac{ e^{D_\varepsilon (x-y)}}{D_\varepsilon -D_0}g(y) dy
+\frac{1}{\varepsilon}O(e^{-(D_\varepsilon -D_0)x})
\int_0^L\frac{ e^{D_\varepsilon (x-y)}}{D_\varepsilon -D_0}g(y) dy\\
& = O(1)\int_x^Le^{D_\varepsilon (x-y)}g(y) dy
+O(1)e^{D_0 x}\int_0^L e^{-D_\varepsilon y}g(y) dy.
\end{split}\]
It is trivial to derive that
\[I_2=O(\varepsilon)e^{-D_\varepsilon(L-x)}
\int_0^LD_0 e^{D_0(L-y)}g(y)dy,\]
and
\[I_4=(1+O(\varepsilon))
\int_0^xe^{D_0(x-y)}g(y)dy.\]
Then the estimates on $I_i$, $i=1, 2, 3, 4$, imply \eqref{comp2}. $\hfill\square$
%----------------------------------------------------------------
\end{appendices}
%----------------------------------------------------------------
%================================================================

%================================================================
%\section*{Acknowledgement}
%
%

\bibliographystyle{apalike} % other choices: plainnat
\bibliography{PancRef}

\end{document}